\documentclass[12pt]{article}

\usepackage{amsmath}
\usepackage[abbrev,msc-links]{amsrefs} 
\usepackage{amssymb}
\usepackage{amsthm}
\usepackage{enumerate}
\usepackage{mathtools}

\usepackage{pgf,tikz,pgfplots}
\usepackage{subcaption}
\usetikzlibrary{arrows}

\renewcommand\proof{\noindent\textsl{Proof. }}
\newcommand\sqr[2]{{\vbox{\hrule height.#2pt
    \hbox{\vrule width.#2pt height#1pt \kern#1pt
        \vrule width.#2pt}\hrule height.#2pt}}}
\renewcommand\qed{%
	\ifmmode\eqno\sqr53
	\else\nolinebreak\ \hfill\sqr53\medbreak\fi}

\title{Oriented or signed Cayley graphs with all eigenvalues 
integer multiples of $\sqrt{\Delta}$}

\date{}
\author{Chris Godsil
\footnote{
C. Godsil gratefully acknowledges the support of the Natural Sciences and Engineering Council of Canada (NSERC), Grant No. RGPIN-9439.}, Xiaohong Zhang\\
University of Waterloo, Waterloo, Ontario, Canada}

\newtheoremstyle{plainsl}%
    {\topsep}
    {\topsep}
    {\slshape} 
    {}
    {\normalfont\bfseries}
    {.}
    { }
    {}

\newcommand{\Z}{\mathbb{Z}}

\DeclareMathOperator{\Gal}{Gal}
\DeclareMathOperator{\lcm}{lcm}
\DeclareMathOperator{\arc}{Arc}
\DeclareMathOperator{\rk}{rank}

\newcommand\Ga{\Gamma}
\newcommand\sg{\sigma}
\newcommand{\cB}{\mathcal {B}}
\newcommand\cA{{\mathcal A}}
\newcommand\cF{{\mathcal F}}
\newcommand\cL{{\mathcal L}}

\newcommand\cx{{\mathbb C}}

\newcommand\rats{{\mathbb Q}}

\newcommand\txtsl[1]{\index{#1}\textsl{#1}}

\newcommand\ip[2]{\langle#1,#2\rangle}
\DeclareMathOperator\elsm{sum}
\DeclareMathOperator{\tr}{tr}

\newcommand\htau{{\hat\tau}}
\newcommand\hsg{{\hat\sigma}}
\newcommand\lca{L[\cA]}

\newcommand\one{{\bf1}}

\newcommand\comp[1]{{\mkern2mu\overline{\mkern-2mu#1}}}
\newcommand\diff{\mathbin{\mkern-1.5mu\setminus\mkern-1.5mu}}

{\theoremstyle{plainsl}
\newtheorem{theorem}{Theorem}[section]
\newtheorem{lemma}[theorem]{Lemma}
\newtheorem{corollary}[theorem]{Corollary}}
{\theoremstyle{remark}
\newtheorem{example}{Example}}
{\theoremstyle{remark}
\newtheorem*{remark}{Remark}}
\newtheorem{question}{Question}

\begin{document}
\maketitle

\section*{Abstract} 

Let $G$ be a finite abelian group. 
Bridges and Mena characterized the Cayley graphs of $G$ that have only integer eigenvalues. 
Here we consider the $(0,1,-1)$ adjacency matrix of an oriented Cayley graph 
or of a signed Cayley graph $X$ on $G$. 
We give a characterization of when all the eigenvalues of $X$ are integer multiples 
of $\sqrt{\Delta}$ for some square-free integer $\Delta$. 
These are exactly the oriented or signed Cayley graphs on which the continuous 
quantum walks are periodic, 
a necessary condition for walks on such graphs to admit perfect state transfer. 
This also has applications in the study of uniform mixing on oriented Cayley graphs, 
as the occurrence of local uniform mixing at vertex $a$ in an oriented graph $X$ 
implies periodicity of the walk at $a$. 
We give examples of oriented Cayley graphs which admit uniform mixing or multiple state transfer. 

\tableofcontents

\section{Introduction}

Formally a \textsl{signed graph} is a graph with a partition of its edges into two
classes (one of which might be empty). We can represent a signed graph by 
a \textsl{signed adjacency matrix}, i.e., a symmetric matrix with non-zero entries $1$ or $-1$. 
An \textsl{oriented graph} is a graph where each edge is assigned a direction,
and can be represented by a skew-symmetric matrix with non-zero entries $1$ and $-1$.
We will refer to this matrix as a \textsl{skew adjacency matrix}.

We are interested in the eigenvalues of signed or skew adjacency matrices of
Cayley digraphs for abelian groups. One reason for our interest is because of connections to quantum walks but, before
addressing this we discuss Cayley graphs and present our main results.

Let $(G,+)$ be an abelian group. The \textsl{Cayley digraph} $X(G,C)$ on $G$ is the digraph with vertex set $G$ 
and with an arc from vertex $g$ to $h$ if and only if $h-g\in C$. 
Let $-C$ denote
\[
	\{-c \, | \, c\in C\}. 
\]
If $C=-C$, then $X(G,C)$ is a Cayley graph (undirected). 
If $C\cap -C=\emptyset$, then $X(G,C)$ is an \textsl{oriented Cayley graph}. 
Otherwise, $X$ is called a \textsl{mixed graph}, with a mixture of edges (corresponding 
to $C\cap (-C)$) and of arcs (corresponding to $C\diff (-C)$). 
Let $C(1)$ and $C(-1)$ be two disjoint inverse-closed subsets of $G$ that do not contain 0, that is, 
\[
0\notin C(1)\cup C(-1),\, 
C(1)=-C(1), \, 
C(-1) =-C(-1),\,
C(1)\cap C(-1)=\emptyset. 
\]
Then the \textsl{signed Cayley graph} $X(G, C(1),C(-1))$ on $G$ is the 
signed graph whose underlying graph is the Cayley graph $X(G, C(1)\cup C(-1))$,  
with edges corresponding to $C(1)$ receiving a positive sign and edges corresponding 
to $C(-1)$ receiving a negative sign. 

We present some of our results. 
We determine when all eigenvalues of an oriented Cayley 
graph $X(G,C)$ or a signed Cayley 
graph $X(G, C(1), C(-1))$ are integer multiples of $\sqrt{\Delta}$ 
for some square-free integer $\Delta$ ($\Delta<0$ for oriented ones and $\Delta>0$ for signed ones). 
More generally, we consider $4$-th root of unity weighted Cayley digraphs, 
which is a combination of an oriented Cayley $X(G,C)$ and a signed 
Cayley $X(G, C(1), C(-1))$, with arcs corresponding to $C$ receiving weight $i$. 
The associated matrix is Hermitian. 

The set of feasible $\Delta$s for a given abelian group $G$ is determined by 
the exponent of $G$, the smallest positive integer $m$ such that $mg=0$ 
for all $g\in G$.  For each feasible $\Delta$, we give a characterization of $C$, 
or of $C(1)$ and $C(-1)$, such that all eigenvalues of the corresponding oriented or signed Cayley graph 
are integer multiples of $\sqrt{\Delta}$. 
This generalizes the classical result of Bridges and Mena \cite{BMinte} on when a Cayley graph 
on an abelian group has only integer eigenvalues. 

For a nonabelian group $G$, there is no known similar characterization of when 
an oriented or signed Cayley graphs for $G$ has all its eigenvalues  integer multiples of $\sqrt{\Delta}$ for some square-free integer $\Delta$. 
But when restricted to normal Cayley graphs for $G$, 
we show that there is a similar (to the abelian case) characterization of the connection set. In the case of nonabelian normal Cayley graphs, the group exponent only provides us with the set of all the possible $\Delta$s, and the feasibility of 
each $\Delta$ is dependent on the group structure (the character values of the group). 

Association schemes play an important role in our analysis. 
In particular, let $\cA$ be the group scheme on an abelian group $G$.
We use automorphisms of the splitting field $L$ of the scheme 
to construct 
automorphisms of the Bose-Mesner algebra $L[\cA]$. 
We make use of the connection between the subalgebra of $L[\cA]$ 
that consists of matrices all of whose entries and eigenvalues 
are in $\rats(\sqrt{\Delta})$ 
and the Galois group $\Gal(L/\rats(\sqrt{\Delta}))$ to define equivalence relations on the group $G$, and make use of Chinese Remainder Theorem to derive
an expression of $\sqrt{\Delta}$ as a linear combination of roots of unity.  

These results have applications to continuous time quantum walks. 
We will see that the oriented (signed) Cayley graphs with the desired eigenvalue 
form are exactly the graphs where the continuous quantum walks are periodic at a 
vertex (in fact periodic at the whole graph here), which is a necessary condition for perfect state transfer to occur in such a graph. 
Furthermore, for (local) uniform mixing to occur on an oriented Cayley graph, 
the graph must be periodic. 
Thus our characterization narrows down the set of 
oriented or signed Cayley graphs on which interesting quantum information transfer phenomenon may occur. 

We give examples of oriented Cayley graphs with uniform mixing that are not constructed 
from Cartesian products of smaller graphs with uniform mixing, 
and examples of multiple state transfer, that is, 
there is a subset $S$ of $V(X)$ such that $|S|>2$ and that between any two vertices 
of $S$ there is perfect state transfer.

\section{Association schemes}

In this section, we review some basics of association scheme. 
If $B$ and $C$ are two $m\times n$ matrices, 
their \textsl{Schur product} $B\circ C$ is the $m\times n$ matrix $B\circ C$ with 
\[
	(B\circ  C)_{jk}=B_{jk}C_{jk}.
\]

An \textsl{association scheme} on $v$ vertices with $d$ classes is a 
set $\cA=\{A_0,\ldots, A_d\}$ of $v\times v$ $(0,1)$-matrices such that 
\begin{enumerate}
\item[(a)] $A_0=I$,
\item[(b)] $\sum_{i=0}^d A_i = J$,
\item[(c)] $A_i^T \in \cA$ for each $i$,
\item[(d)] $A_iA_j=A_jA_i $ for all $i,j=0,\ldots, d$.
\item[(e)] $A_iA_j\in \operatorname{span}(\cA)$ for all $i,j=0,\ldots, d.$
\end{enumerate}
For example, if $X$ is a strongly regular graph with 
complement $\comp{X}$, then the matrices
\[
	\{A_0=I, A_1=A(X), A_2=A(\bar{X})\}
\] 
form an association scheme with two classes. 
We will make use of a special type of association scheme, 
the group scheme of a finite abelian group. 
See Section~\ref{sec:BMena} for details. 

Since $J$ and $A_i$ commute, $A_i$ has constant row sums and constant column sums, and the two sums are equal, which we denote by $v_i$. 

Let $\cx[\cA]$ be the vector space over $\cx$ spanned by $\cA$. 
It is closed under both the matrix multiplication and Schur multiplication,  
and is called the \textsl{Bose-Mesner algebra}  of $\cA$ \cite{Godsilcombi}.  
On $\cx[\cA]$, or more generally, 
on $M_{v\times v}(\cx)$, 
the set of all $v\times v$ complex matrices, 
\[
\langle M_1, M_2 \rangle = \tr(M_1^*M_2)
\]
is an inner product. 
Let $M$ be a complex matrix, 
we use $\comp{M}$ to denote the complex conjugate of $M$, 
and use $\elsm (M)$ to denote the sum of all entries of $M$. 
Then 
\begin{equation}\label{eq:innerpro}
\langle M_1, M_2 \rangle = \tr(M_1^*M_2)
	=\elsm(\overline{M_1}\circ M_2).
\end{equation}
Since $A_i\circ A_j=\delta_{ij}A_i$, 
we have $\langle A_i, A_j \rangle = \delta_{ij} \elsm(A_i)=\delta_{ij} v v_i$, 
where as above, $v_i$ is the constant row sum of $A_i$. 
Therefore $\cA$ is an orthogonal basis of $\cx[\cA]$ consisting of pairwise orthogonal Schur idempotents. 
In fact, there is another important orthogonal basis of $\cx[\cA]$.

\begin{theorem}\label{thm:prinidem}\cite{Godsilcombi}
Let $\cA=\{A_0,\ldots, A_d\}$ be an association scheme with 
$d$ classes on  $v$ vertices. 
Then there is a set of pairwise orthogonal idempotent matrices $E_0,\ldots, E_d$ and numbers $p_i(j)$ such that 
\begin{enumerate}
\item[(a)] $\sum_{j=0}^dE_j=I$,
\item[(b)] $A_iE_j=p_i(j)E_j$ for all $i,j=0,\ldots, d$,
\item[(c)] $E_0=\frac1v J$, $E_j^*=E_j$ for $j=0,\ldots, d$, 
\item[(d)] $E_0,\ldots, E_d$ is a basis for $\cx[\cA]$.\qed
\end{enumerate}
\end{theorem}
$E_0,\ldots, E_d$ are called the \textsl{principal idempotents} of the scheme. 
The eigenvalues $p_i(j)$, $j, i=0,\ldots, d$ are called the \textsl{eigenvalues of the association scheme}. 
More details about association schemes, including the following parameters, can be found in Chapter 12 of \cite{Godsilcombi}. 

Now we introduce two more families of parameters associated to an association scheme. 
Since $\cx[\cA]$ is closed under Schur multiplication and has $E_0,\ldots, E_d$ as a basis, 
there are numbers $q_{ij}(k)$ such that 
\[
E_i \circ E_j =\frac1v \sum_{k=0}^d q_{ij}(k)E_k.
\]
The numbers $q_{ij}(k)$, $i,j,k=0,\ldots, d$, are called \textsl{Krein parameters} of the association scheme. 
Note that these parameters can be written in terms of the eigenvalues 
\begin{equation}\label{eq:kreineigen} 
q_{ij}(k)=vm_im_j\sum_{r=0}^d p_r(i)p_r(j)p_r(k)/\big(\elsm(A_r)\big)^2,
\end{equation}
where $m_i=\rk(E_i)$. 
Making use of the fact the matrices $A_0,\ldots, A_d$ form a basis of $\cx[\cA]$, 
we know that there exist numbers $q_i(j)$ such that 
\[
E_i=\frac1v \sum_{j=0}^d q_i(j) A_j.
\]
The values $q_i(j)$ for $j,i=0,\ldots, d$ are called the \emph{dual eigenvalues} of the scheme. 
Eigenvalues and dual eigenvalues can be written in terms of each other as 
\begin{equation} \label{eq:eigendualeig}
p_i(j)m_j=q_j(i)v_i.  
\end{equation}

We have seen that, given an association scheme, there is a Bose-Mesner algebra associated to it. 
On the other hand, given a matrix algebra, we can check whether it is the Bose-Mesner algebra of some association scheme. 
\begin{lemma}\label{lem:coheBM}\cite{comcoheBS}
Let $\mathcal{W}$ be a commutative matrix algebra of order $n$ over a subfield $E$ of $\cx$. 
If   
\begin{enumerate}
\item[(a)] $I,J\in  \mathcal{W}$,
\item[(b)] If $M\in  \mathcal{W}$, then $M^T\in  \mathcal{W}$,
\item[(c)] $MN\in  \mathcal{W}$ and $M\circ N\in  \mathcal{W}$ for all $M,N\in  \mathcal{W}$,
\end{enumerate}
then there exists an association scheme $\cA$ such that $\mathcal{W}=E[\cA]$, 
the span of $\cA$ over $E$.\qed
\end{lemma}

\section{Bose-Mesner algebra automorphisms}
\label{BMauts}

Let $\cA=\{A_0,\ldots, A_d\}$ be an association scheme. 
The \textsl{splitting field} of
$\cA$ is the extension $L$ of the rational field $\rats$ generated by the
eigenvalues of the scheme.  From equation \eqref{eq:eigendualeig}, a relation
between the dual eigenvalues and the eigenvalues, we see that the
splitting field is also generated by the dual eigenvalues.  

We have introduced the Bose-Mesner algebra over the complex field $\cx$, 
in fact, such an algebra can be defined over a subfield of $\cx$. 
Let $E$ be any extension field of $L$ that is closed under complex conjugation, 
and let $E[\cA]$ be the vector space over $E$ spanned by $\cA$. 
Then $E[\cA]$  is closed under both  matrix multiplication, Schur multiplication,  
and complex conjugation, 
and we call it the \textsl{Bose-Mesner algebra} of $\cA$ over $E$. 
In this section, 
we define Bose-Mesner algebra automorphisms of $E[\cA]$ and derive some related results. 
In particular, 
we show that for any Bose-Mesner algebra automorphism $\psi$ of $E[\cA]$, 
the set of fixed-points of $E[\cA]$ under $\psi$ is the Bose-Mesner algebra of a subscheme of $\cA$. Some of the results generalize results in \cite{godsil2010generalized}. 

Note that most of the definitions and results in this section can be generalized 
to a  more general field, 
but in this work we focus on subfield $E$ of $\cx$ that contains all the eigenvalues 
of the scheme and is closed under complex conjugation.  
In this case, $E[\cA]$ shares most properties of $\cx[\cA]$: 
we can talk about both the Schur idempotents and principal idempotents in the algebra, 
we can take complex conjugate and conjugate transpose of matrices, and 
Equation~\eqref{eq:innerpro} provides an inner product on 
the algebra. 

Let $\cA$ be an association scheme with splitting field $L$. 
Let $E$ be an extension of $L$ that is closed under complex conjugation. 
An $E$-linear map $M\mapsto M^\psi$ on $E[\cA]$ 
is a \txtsl{Bose-Mesner algebra automorphism} if for all $M$ and $N$ 
in $E[\cA]$:
\begin{flalign}\label{eq:BSauto}
& \text{(a)} \;\; (MN)^\psi =M^\psi N^\psi,& \nonumber \\ 
&\text{(b)} \;\; (M\circ N)^\psi =M^\psi\circ N^\psi, &  \\
&\text{(c)}  \;\; \psi \text{ is invertible}, & \nonumber \\
&\text{(d)} \;\;  (M^*)^\psi =(M^\psi)^*. \nonumber
\end{flalign}

Now we show that if a $\rats$-linear map of $E[\cA]$ satisfies both (a) and (b), 
then it also satisfies (c).

\begin{lemma}\label{lem:qlinproshurpro}
Let $\cA$ be an association scheme. 
Denote the splitting field of $\cA$ by $L$ and let $E$ be an extension field of $L$. 
Let $\psi$ be a non-zero $\rats$-linear mapping on $E[\cA]$ 
that preserves matrix multiplication and Schur multiplication. 
Then $\psi$ is invertible. 
In particular, 
it permutes elements of $\cA$ and permutes the principal idempotents $E_j$s, 
with $I^\psi = I$ and $J^\psi =J$. 
\end{lemma}

\proof 
Since $\psi$ preserves the two types of matrix multiplications, 
it follows immediately that $\psi$ maps Schur idempotents to Schur
idempotents and matrix idempotents to matrix idempotents.  
Hence $J^\psi$ is a 01-matrix.   
Now we prove that $\psi$ is invertible. 
Assume $\cA$ is on $v$ vertices. Since $J^2=vJ$ and $\psi$ is $\rats$-linear and nonzero,  we have
$$
(J^\psi)^2=(J^2)^\psi=vJ^\psi;
$$
it follows that $J^\psi=J$.  Consequently
\begin{equation}\label{eq:permai}
J=J^\psi=\sum_i A_i^\psi.
\end{equation}
Since $JA_i=A_iJ$ we know that $A_iJ=v_iJ$ for some positive integer $v_i$, and 
\[
A_i^\psi J =(A_iJ)^\psi =v_i J. 
\]
Therefore  $A_i^\psi\neq 0$ for any $i=0,\ldots, d$;   
this combined with \eqref{eq:permai} and the fact $A_i^\psi$ is a Schur idempotent  implies that $\psi$ permutes the elements of $\cA$.
Therefore $\psi$ maps a basis of $E[\cA]$ to a basis of $E[\cA]$ and 
 it is invertible. 
Now from 
\[
E_i^\psi E_j^\psi = (E_iE_j)^\psi=\delta_{ij}E_i^\psi,
\]
we know $\{E_0^\psi,\ldots, E_d^\psi\}$ is a set of mutually orthogonal idempotents. 
Since $\psi$ is invertible, 
none of the $E_j^\psi$ is 0, 
and therefore $\psi$ permute the idempotents $E_j$s. 
From this we get 
\[
	I^\psi = (\sum_j E_j)^\psi = \sum_j E_j^\psi = I.\qed
\]

If further the above $\psi$ is $E$-linear, 
then in addition to being invertible, it also satisfies (d) in \eqref{eq:BSauto}:
\begin{lemma}\label{lem:clinproschurpro*}
Let $\cA$ be an association scheme with splitting field $L$ 
and let $E$ be an extension field of $L$ that is closed under complex conjugation. 
If $\psi$ is a non-zero $E$-linear mapping on $E[\cA]$ 
that preserves matrix multiplication and Schur multiplication, 
then $(M^T)^\psi =(M^\psi)^T$ and $(M^*)^\psi =(M^\psi)^*$. 
\end{lemma}
\proof  
Since $E$ is closed under complex conjugation and $\cA$ is closed under transposition, 
we know that if $M\in E[\cA]$, then $M^T,\overline{M}, M^*\in E[\cA]$. 
Therefore \eqref{eq:innerpro} is an inner product on $E[\cA]$, 
with  
both $\cA$ and $\{E_0,\ldots, E_d\}$ are orthogonal bases. 
Since
\begin{equation}\label{eq:inprosch}
\elsm(I \circ (A_iA_j)) =\tr(A_iA_j)=\ip{A_i^T}{A_j},
\end{equation}
we find that $\elsm (I \circ (A_iA_j))\ne0$ if and only if $A_j=A_i^T$.  
As earlier, denote the constant row sum of $A_k$ as $v_k$. 
Since $\psi$ is $E$-linear and preserves matrix multiplication and Schur 
multiplication,  Lemma~\ref{lem:qlinproshurpro} implies $I^\psi =I$ and 
that $A_k^\psi, (A_k^T)^\psi\in \cA$ for any $k=0,\ldots, d$. 
Now  
\[
\elsm \Big(I \circ \big(A_k^\psi (A_k^T)^\psi \big)\Big) =
\elsm \Big(\big(I \circ (A_kA_k^T)\big)^\psi \Big)=
\elsm(v_k I)=
v_kv\neq0.  
\]
From the characterization of when \eqref{eq:inprosch} is nonzero with 
$A_i=A_k^\psi$ and $A_j=(A_k^T)^\psi$, 
we know $(A_k^\psi)^T= (A_k^T)^\psi$ for any $k$. 
From the linearity of $\psi$, 
we have $(M^T)^\psi =(M^\psi)^T$ and $(M^*)^\psi =(M^\psi)^*$ for any $M\in E[\cA]$. \qed

Combining the above results, 
we have 
\begin{theorem}\label{thm:BMauto}
Let $\cA$ be an association scheme with splitting field $L$ and let $E$ be an 
extension field of $L$ that is closed under complex conjugation. 
A non-zero $E$-linear map 
$M\mapsto M^\psi$ is a Bose-Mesner algebra automorphism of $E[\cA]$  
if and only if for all $M$ and $N$ 
in $E[\cA]$:
\begin{enumerate}
\item[(a)] $(MN)^\psi =M^\psi N^\psi$, and 
\item[(b)] $(M\circ N)^\psi =M^\psi\circ N^\psi$. \qed
\end{enumerate}
\end{theorem}

The transpose map is a Bose-Mesner algebra automorphism on $E[\cA]$, which is non-trivial if
$\cA$ is not symmetric, that is, if not all matrices in $\cA$ are symmetric.

Let $\pi=\{C_0=\{0\},C_1,\ldots, C_e\}$ be a partition of the indices $\{0,\ldots, d\}$, 
and let 
\[
B_i=\sum_{j\in C_i}A_j,  \quad i=0,\ldots, e.
\]
That is, $\cB=\{B_0,\ldots, B_e\}$ is obtained by merging classes of $\cA$. 
If $\cB$ is an association scheme, 
then it is called a \textsl{subscheme} (or  \textsl{fusion scheme}) of $\cA$.  
We are going to use Bose-Mesner algebra automorphisms to construct subschemes.  If
$\psi$ is a Bose-Mesner algebra automorphism on $E[\cA]$, the \textsl{fixed-point space} of
$\psi$ is the set of matrices in $E[\cA]$ that are fixed by $\psi$.
This is evidently a subspace of $E[\cA]$. Further, 

\begin{lemma}
Let $\cA$ be an association scheme with splitting field $L$.
Let $E$ be an extension of $L$ that is closed under complex conjugation. 
The fixed-point space of a Bose-Mesner algebra automorphism of $E[\cA]$ 
is the Bose-Mesner algebra of a subcheme of $\cA$. 
\end{lemma}

\proof
The fixed-point space of a Bose-Mesner algebra automorphism is closed under multiplication, Schur
multiplication, transpose and contains $I$ and $J$. 
It is also commutative with respect to matrix multiplication, as $E[\cA]$ is. 
By Lemma~\ref{lem:coheBM}, 
it is the Bose-Mesner algebra of a scheme $\cB$. 
Since for any $B_i\in \cB$, $B_i\in E[\cA]$, 
we know $B_i$ is a sum of matrices in $\cA$ and $\cB$ is a subscheme of $\cA$.\qed

By way of example, consider the transpose map acting on $\cA$.  Its
fixed-point space is spanned by those Schur idempotents that are
symmetric, together with the matrices
$$
A_i+A_i^T,
$$
where $A_i$ is not symmetric.  By the lemma, these matrices are the
Schur idempotents of a symmetric subscheme of $\cA$.

\section{Galois}\label{mungal}

Let $\cA$ be an association scheme with splitting field $L$. 
Note that $L$ is closed under complex conjugation: 
if $A_iE_j=p_i(j)E_j$, then $A_i^TE_j=\overline{p_i(j)}E_j$. 
From now on, we focus on the  Bose-Mesner algebra $L[\cA]$.  
The \textsl{Krein field} $K$ is the extension
of $\rats$ generated by the Krein parameters.  
It follows from equation \eqref{eq:kreineigen} 
that the Krein field $K$ is a subfield of $L$.  
In this section, we construct 
two group actions of the Galois group $\Gal(L/\rats)$ over $\lca$.  
By making use of $K$, we give a characterization of 
when $\tau\in\Gal(L/\rats)$  is a Bose-Mesner algebra automorphism 
of $\lca$ through the second action. 

Let $\sg$ be a field automorphism of $L$ and let $\alpha\in L$, the the image of $\alpha$ under $\sg$ is denoted by $\alpha^\sg$ or $\sg(\alpha)$.

\medbreak

Now we define the first action of $\Ga=\Gal(L/\rats)$ over $\lca$. 
Let $\sg\in\Ga$ and $M\in\lca$, 
define $M^\sg$ by 
\begin{equation*}
(M^\sg)_{jk}=(M_{jk})^\sg. 
\end{equation*}
That is, 
$M^\sg$ is the matrix obtained by
applying $\sg$ to each entry of $M$.  
This gives the \textsl{entrywise} action of $\Ga$.  
Note that for any $M,N\in \lca$, 
\[
(MN)^\sg=M^\sg N^\sg,\quad (M\circ N)^\sg=M^\sg \circ N^\sg. 
\]
Further, this action is linear over $\rats$, but not over $L$. 
By Lemma~\ref{lem:qlinproshurpro}, 
this entrywise action permutes the elements of $\cA$ (it in fact fixed all the Schur idempotents) and 
permutes the principal idempotents of $\cA$.

We define a second action of $\Ga$ on $\lca$.  Suppose $\sg\in\Ga$
and $M\in\lca$.  
Then $M=\sum_j a_j E_j$ for some $a_j\in L$. 
Define $M^\hsg$ by
\begin{equation}\label{eq:fieldlinact}
M^\hsg= \sum_j a_j E_j^\sg,
\end{equation}
where $ E_j^\sg$ is the image of $E_j$ under the above entrywise action of $\sg$ (note that $E_j^\hsg=E_j^\sg$). 
This is an $L$-linear map on $\lca$. 
In fact, it also preserves matrix multiplication: 

\begin{lemma}\label{lem:hatpremulti}
Let $\cA$ be an association scheme and $L$ be its splitting field. 
For $\sg\in\Gal(L/\rats)$, let $\hsg$ be the $L$-linear mapping on $L[\cA]$ defined as above in \eqref{eq:fieldlinact}. 
Then for any $M,N\in\lca$,  
$(MN)^\hsg=M^\hsg N^\hsg$. 
\end{lemma}

\proof
Assume $\cA=\{A_0,\ldots, A_d\}$ and that $E_0,\ldots, E_d$ are its 
principal idempotents. 
As discussed above, the entrywise action of $\sg$ on $L[\cA]$ preserves 
matrix multiplication and Schur multiplication, 
and satisfies  
$\{E_0^\sg,\ldots, E_d^\sg\}=\{E_0,\ldots, E_d\}$ 
and $A_i^\sg=A_i$ for $i=0,1,\ldots, d$. 
Now 
\[
(E_iE_j)^\hsg=\delta_{i,j}E_i^\sg=E_i^\sg E_j^\sg
=E_i^\hsg E_j^\hsg, 
\]
by $L$-linearity of $\hsg$, we have for any $M,N\in \lca$,
\[
	(MN)^\hsg= M^\hsg N^\hsg.\qed
\]

Now we check when $\hsg$ preserves Schur multiplication. 
\begin{theorem}\label{eigaut}
Let $\cA$ be an association scheme with splitting field $L$ and Krein
field $K$.  
Let  $\sg\in \Gal(L/\rats)$,
then $\hsg$, as defined in \eqref{eq:fieldlinact}, is a Bose-Mesner algebra 
automorphism of $\lca$ if and only if $\sg$ fixes each element of $K$ (that is, $\sg\in \Gal(L/K)$). 
\end{theorem}

\proof
Since $\hsg$ preserves matrix multiplication (Lemma~\ref{lem:hatpremulti}) and is  $L$-linear, 
by Theorem~\ref{thm:BMauto}, 
it suffices to prove that $\hsg$  
preverses Schur multiplication if and only if $\sg$ fixes all the Krein parameters. 
Take any $M,N\in\lca$.  On the one hand, 
\[
(E_i\circ E_j)^\hsg =\left(\frac{1}{ v}\sum_r q_{i,j}(r) E_r\right)^{\!\! \hsg}
    =\frac{1}{ v}\sum_r q_{i,j}(r) E_r^\sg. 
\]
While, on the other, since the entrywise action of $\sg$ preserves 
Schur multiplication and is $\rats$-linear,   
\[
E_i^\hsg\circ E_j^\hsg =E_i^\sg\circ E_j^\sg =(E_i\circ E_j)^\sg
    =\frac{1}{ v}\sum_r q_{i,j}(r)^\sg E_r^\sg.
\]
Comparing these two equations yields that
\[
(E_i\circ E_j)^\hsg =E_i^\hsg\circ E_j^\hsg
\]
for all $i$ and $j$, if and only if $\sg$ fixes each Krein parameter.
\qed

Using related, but distinct, actions of $\Gal(L/K)$,
Munemasa \cite{Mune} proved that $\Gal(L/K)$ lies in the centre of $\Gal(L/\rats)$.
(Similar results appear in \cite{deBoeGoeree,CosteGannon}.)  
Since the argument is short, we present a version of it here.  
Take any $\sg\in\Gal(L/\rats)$ and $\tau\in \Gal(L/K)$. 
Then by Theorem~\ref{eigaut}, 
$\htau$ is a Bose-Mesner algebra automorphism on $L[\cA]$ and therefore permutes the $A_j$s. 
Since the dual eigenvalues also generate $L$, 
to prove $\sg\tau=\tau\sg$ on the splitting field $L$, 
we just need to show that $\sg\tau=\tau\sg$ on all the  dual eigenvalues of the scheme. 
Since $\sg$ fixes the $A_i$s, $\htau$ permutes the $A_i$s and is $L$-linear, it follows that 
\begin{align*}
(vE_i)^{\sg\htau}
& =\left(\sum_j q_i(j)A_j\right)^{\!\!\sg\htau}
=\left(\sum_j q_i(j)^\sg A_j\right)^{\!\!\htau}
=\sum_j q_i(j)^\sg A_j^\htau \\
&= \left(\sum_j q_i(j) A_j^\htau\right)^{\!\!\sg}
    =\left(\sum_j q_i(j)A_j\right)^{\!\!\htau\sg}
    =(vE_i)^{\htau\sg}.
\end{align*}
Therefore $E_i^{\sg\htau}=E_i^{\htau\sg}$ for all $i$. 
Now since $\sg$ permutes the $E_i$s and fixes the $A_i$s, 
\[
E_i^{\sg\htau}=E_i^{\sg\tau}=\frac{1}{ v}\sum_j q_i(j)^{\sg\tau}A_j,
\]
and 
\[
E_i^{\htau\sg}=E_i^{\tau\sg}=\frac{1}{ v}\sum_j q_i(j)^{\tau\sg} A_j.
\]
Combining the last three equations, we conclude that
\[
q_i(j)^{\sg\tau}=q_i(j)^{\tau\sg}
\]
for all $i,j$. 
Therefore $\Gal(L/K)$ lies in the centre of $\Gal(L/\rats)$.

Let $F$ be a subfield of $L$ that contains $K$. 
Consider the subset 
 \[
 \cF=\{M\in L[\cA]\,|\, \text{ all the entries and eigenvalues of $M$ are in $F$}\} 
 \]
 of $L[\cA]$. 
Hou showed that $\cF=F[\cB]$ for some subscheme $\cB$ of $\cA$ 
and gave a characterization of $\cB$ \cite{Hou1992}.  
This result plays an important role in our later analysis. 
Here we provide a different 
 proof using the Bose-Mesner algebra automorphism $\htau$ on $L[\cA]$ for $\tau\in\Gal(L/K)$.

\begin{theorem}\cite{Hou1992}\label{fldsub}
Let $\cA$ be an association scheme with splitting field $L$
and Krein field $K$. 
Consider the action of $\Gal(L/K)$ on $L[\cA]$ defined in \eqref{eq:fieldlinact}. 
Let $F$ be a subfield of $L$ that contains $K$, 
and let $\cF$ be the set of matrices in $\lca$ with
eigenvalues and entries in $F$.
Then $\cF=F[\cB]$, 
where $\cB$ is the subscheme of $\cA$ fixed by 
$\{\htau\, | \, \tau\in \Gal(L/F)\}$.
\end{theorem}

\proof
Clearly $\cF$ is a transpose-closed algebra over $F$ (under the usual matrix multiplication). 
Take $M,N\in \cF$ and assume 
\[
M=\sum_i a_i E_i,\qquad N=\sum_i b_i E_i.
\]
Then $a_i,b_i\in F$ and entries of $M\circ N$ are in $F$. Further   
\[
M\circ N= \sum_{i,j} a_i b_j E_i\circ E_j=\sum_{k}(\frac1v\sum_{i,j}a_ib_jq_{ij}(k)) E_k 
\]
and the fact $q_{ij}(k)\in K\subseteq F$ imply all 
the eigenvalues of $M \circ N$ lie in $F$. 
Therefore $\cF$ is Schur-closed.  
Note that $I,J\in \cF$. 
By Lemma~\ref{lem:coheBM}, 
$\cF=F[\cB]$ for some subscheme $\cB$ of $\cA$. 
Take any $\tau\in \Gal(L/K)$ 
and  any $M=\sum_i a_i E_i \in\lca$, 
we have 
\[
M^{\tau^{-1}\htau}=\sum_i a_i^{\tau^{-1}}E_i.
\]
If $M$ is a rational matrix, 
then $M^{\tau^{-1}}=M$. This shows that a rational matrix in $\lca$ is fixed by $\htau$ if and
only if all its eigenvalues are fixed by $\tau$.   
Since a rational matrix
lies in $\cF$ if and only if all its eigenvalues are fixed by $\Gal(L/F)$, 
it follows a rational matrix
lies in $\cF$ if and only if it is fixed by $\{\htau\, | \, \tau\in \Gal(L/F)\}$.\qed

\section{Subschemes of a group scheme}
\label{sec:BMena}

In this section we focus on the group scheme $\cA$ 
of a finite abelian group $G$. 
Let $F$ be a subfield of the splitting field $L$ of $\cA$, 
we find the subschemes of $\cA$ whose eigenvalues are all in $F$. 
In particular, we give a characterization on the connection set $C$ for the Cayley digraph
$X(G,C)$ to have all its eigenvalues in $F$, which provides an alternate proof of 
 a classical result by Bridges and Mena on when a Cayley graph is integral.

Assume $|G|=n$. 
For any $g\in G$, 
let $A_g$ be the $n\times n$ (0,1) matrix with rows and columns indexed by elements of $G$, 
such that $(A_g)_{u,v}=1$ if and only if $v-u=g$. 
So $A_g$ is a permutation matrix and $A_gA_h=A_{g+h}=A_hA_g$. 
In fact, 
\[
\cA = \{A_g\, | \, g\in G\}
\]
is an association scheme on $n$ vertices with $n-1$ classes, 
and it is called the \textsl{group scheme} of $G$.

First consider the cyclic group $G=\Z_n$. 
Let $P_n$ be the $n\times n$ permutation matrix 
\[
P_n=\begin{bmatrix} 
0 & 1 & 0 & \cdots & 0\\
0 & 0 & 1 &\cdots & 0\\
0 & \vdots & \vdots & \ddots& \vdots\\
0 & 0 & 0 & \cdots & 1\\
1 & 0 & 0 &\cdots  &0
\end{bmatrix}. 
\]
Then 
\[ 
\mathcal{A}=\{A_{0}=I, A_{1} = P_n, \ldots, A_{n-1}=P_n^{n-1}\}
\]
form an association scheme on $n$ vertices with $n-1$ classes,
and is called a cyclic group scheme. 
Let $\zeta_n=e^{i\frac{2\pi}{n}}$. 
For $r=0,1,\ldots, n-1$,  
let $E_r$ be the $n\times n$ matrix with 
$(E_r)_{j,k}=\frac1n \zeta_n^{r(j-k)}$ for $j,k=0,\ldots, n-1$. 
Then $E_0,\ldots, E_{n-1}$ are the principal idempotents  of $\cA$. 
Since the eigenvalues of the scheme are $p_i(j)=\zeta_n^{ij}$, 
the splitting field $L$ is the cyclotomic field obtained by adding the 
primitive $n$-th root $\zeta_n$ to $\rats$, which we denote by $\rats(\zeta_{n})$. 
Now  $\Ga = \Gal(L / \rats) \cong \Z_n^*$, the multiplicative group of modulo $n$,
where $k\in \Z_n^*$ corresponds to $\tau_k\in \Gal(L / \rats)$ such that 
\begin{equation*}
\tau_k(\zeta_n)=\zeta_n^k.
\end{equation*}
Denote this group isomorphism from $\Z_n^*$ to $\Gal(L / \rats)$ by $\phi$. 

Since $E_i\circ E_j=\frac1n E_{i+j}$, 
where the summation of subscripts is over $\Z_n$, 
all the Krein parameters of $\cA$ are rational, 
and therefore the Krein field satisfies $K=\rats$. 
By Theorem~\ref{eigaut}, for any $\tau_k\in \Gal(L/K)=\Gal(L/\rats)\cong \Z_n^*$, 
$\htau_k$ is a Bose-Mesner algebra automorphism on $L[\cA]$. 
Further, 
it can be checked for any $k\in \Z_n^*$,  
\[
E_m^{\hat{\tau}_k}=E_{mk} \text{ and } 
A_m^{\hat{\tau}_k} = A_{mk^{-1}}, 
\]
where the product and inverse of the subscripts are taken in $\Z_n^*$.

More generally, 
let 
$G = \Z_{n_1}\times \Z_{n_2}\times\cdots \times \Z_{n_e}$,  
a direct product of cyclic groups and let 
\begin{equation}\label{eq:factorscheme}
\cA_1=\{A_{1,0},\ldots, A_{1,n_1-1}\},\ldots, 
\cA_e=\{A_{e,0},\ldots, A_{e,n_e-1}\}
\end{equation}
 be the cyclic group schemes corresponding to 
the cyclic groups $\Z_{n_1},\ldots, \Z_{n_e}$, respectively, 
then the group scheme on $G$ is 
\begin{equation}\label{eq:schemeG}
\cA =\cA_1\times \cdots \times \cA_e=
\{A_{1,j_1}\otimes A_{2,j_2}\otimes\cdots\otimes A_{e,j_e}\; | \, A_{k,j_k} \in \cA_k \}.
\end{equation}
We denote the matrix $A_{1,j_1}\otimes A_{2,j_2}\otimes\cdots\otimes A_{e,j_e}$ 
by $A_{(j_1,\ldots, j_e)}$ for short. 
Further, for $k=1,\ldots, e$,  
let $E_{k,j_k}$, $j_k=0,1,\ldots, n_k-1$, 
be the principal idempotents of the group scheme $\cA_k$. 
Then 
\[
E_{(j_1,\ldots, j_e)}:=E_{1,j_1}\otimes E_{2,j_2}\otimes\cdots\otimes E_{e,j_e}, \; 
j_k=0,\ldots, n_k-1
\]
are the principal idempotents of $\cA$, 
and 
\[
A_{(j_1,\ldots, j_e)} E_{(k_1,\ldots, k_e)} = 
\zeta_{n_1}^{j_ik_1} \cdots \zeta_{n_e}^{j_ek_e}E_{(k_1,\ldots, k_e)}. 
\]
Therefore the splitting field of $\cA$ is $L=\rats(\zeta_{\lcm(n_1,\ldots, n_e)})$, 
and the Krein field is $K=\rats$.

Now we characterize when the adjacency matrix of a Cayley digraph for $G$ 
has all its eigenvalues belonging to a given subfield of $L$. 
When talking about the adjacency matrix of a Cayley digraph $X(G,C)$, 
we mean the $(0,1)$ matrix $A_C$ such that $(A_C)_{g,h}=1$ if and only if $h-g\in C$.  
That is, $A_C=\sum_{g\in C}A_g$.   
Recall that the \textsl{exponent} of an abelian group $G$ is the least 
positive integer $m$ such that $mg=0$ for all $g\in G$.

\begin{theorem} \label{thm:HFconnectionset}
Let $G$ be a finite abelian group with exponent $n_1$. 
Let $\cA$ be the group scheme on $G$. 
Let $F$ be a subfield of $L=\rats(\zeta_{n_1})$. 
Let $\tau_k\in \Gal(L/ \rats)$ be defined as 
$\tau_k(\zeta_{n_1})=\zeta_{n_1}^k$. 
Let $\phi$ be the group isomorphism from $Z_{n_1}^*$ to $\Gal(L/ \rats)$ such that 
$\phi(k)=\tau_k$ and let $Z_F=\phi^{-1}(\Gal(L/ F))$. 
Define an equivalence relation $\sim_F$ on $G$ by $g\sim_F h$ if and only if $g=kh$ for some $k\in Z_F$. 
Then the set $\cF$ of matrices in $L[\cA]$ with entries and eigenvalues in $F$  satisfies $\cF=F[\cB]$, 
where $\cB$ is the subscheme of $\cA$, 
with each element $B_i$ corresponding to an equivalence class of $G$ under $\sim_F$: $B_{[g]_F}=\sum_{h:h\sim_F g}A_h$. 
 In particular,  a Cayley digraph $X(G,C)$  has all its eigenvalues lying in $F$ 
if and only if the connection set $C$ is a union of equivalence classes of $\sim_F$. 
\end{theorem}

\proof
By the fundamental theorem of finite abelian groups, 
we can assume without loss of generality that 
$G = \Z_{n_1}\times \Z_{n_2}\times\cdots \times \Z_{n_e}$ with $n_i\mid n_{i-1}$ for $i=2,\ldots, e$. 
By the above argument, 
we know $\cA$ has splitting field $L=\rats(\zeta_{n_1})$ and Krein field $K=\rats$. 
Since $L\supseteq F\supseteq K= \rats$, by Theorem~\ref{fldsub},  
$\cF=F[\cB]$ for the subscheme $\cB$ of $\cA$ fixed by $\{\htau\, | \, \tau \in \Gal(L/F)\}$. 
Take any $B_i\in \cB$, then $B_i=A_C=\sum_{g\in C}A_g$ for some subset $C\subseteq G$. More precisely, using the notation as in \eqref{eq:schemeG}, 
\[A_C
=\sum_{(j_1,\ldots, j_e)\in C} 
A_{1,j_1}\otimes A_{2,j_2}\otimes\cdots\otimes A_{e,j_e}=\sum_{(j_1,\ldots, j_e)\in C}A_{(j_1,\ldots, j_e)}. 
\]
Now for any $k\in \Z_{n_1}^*$, 
\begin{equation*} 
(A_{(j_1,\ldots, j_e)})^{\hat{\tau}_k}
=A_{k^{-1}(j_1,\ldots, j_e)}\quad \text{ and }\quad 
(E_{(j_1,\ldots, j_e)})^{\hat{\tau}_k}
=E_{k(j_1,\ldots, j_e)}.
\end{equation*}
Therefore $A_C$ is fixed by $\htau_k$ if and only if $C$ is closed under multiplication 
by $k^{-1}\in \Z_{n_1}^*$, 
and $A_C$ is fixed by $\{\htau_k\,|\, \tau_k\in \Gal(L/ F)\}$ if and only if $C$ 
is closed under multiplication by $Z_F$. 
Therefore each Schur idempotent $B_i$ of $\cB$ corresponds to an equivalence 
class of $G$ under $\sim_F$.\qed

Let $G$ be an abelian group. 
Define a relation $\sim$ on $G$ by $g\sim h$ if and only if they generate the same (cyclic)  
subgroup of $G$, that is, $\langle g\rangle = \langle h\rangle$. 
This is an equivalence relation on $G$. 
Denote by $A_{[g]}$ the matrix $\sum_{h:h\sim g}A_h$. 
A direct consequence of the above theorem is the 
 following result by Bridges and Mena which characterizes when a Cayley 
 graph for $G$ has only integer eigenvalues.
\begin{corollary}\cite{BMinte} \label{thm:integtran}
Let $G$ be an abelian group. 
Define an equivalence relation $\sim$ on $G$ by $g\sim h$ if and only if $\langle g\rangle = \langle h\rangle$. 
Then the eigenvalues of the Cayley graph $X(G,C)$ are all 
integers if and only if $C$ is a union of the equivalence classes of $\sim$. 
\end{corollary}

\proof 
Again without loss of generality, 
we assume that $G = \Z_{n_1}\times \Z_{n_2}\times\cdots \times \Z_{n_e}$ 
with $n_i\mid n_{i-1}$ for $i=2,\ldots, e$. 
Therefore the splitting field is $L=\rats(\zeta_{n_1})$ and the Krein field is $K=\rats$. 
Since any rational eigenvalue of an integer matrix is an integer, 
 finding Cayley graphs on $G$ with only integer eigenvalues is equivalent to 
 finding Cayley graphs with only rational eigenvalues. 
With the notion of  Theorem~\ref{thm:HFconnectionset}, 
we are taking the subfield $F$ of $L=\rats(\zeta_{n_1})$ to be $F=\rats$. 
Therefore  $Z_F=\Z_{n_1}^*$, 
and $X(G,C)$ has only integer eigenvalues if and only if $C$ is closed under 
multiplication by $Z_F=\Z_{n_1}^*$. 
Now the result follows from the fact that two group elements $g,h\in G$ generate 
the same cyclic subgroup if and only if $g=kh$ for some $k\in \Z_{n_1}^*$.\qed

\section{Quadratic subfields of $\rats(\zeta_{n})$}

Let $G$ be a finite abelian group and $\cA$ be the associated group scheme. 
Corollary~\ref{thm:integtran} characterizes matrices in $L[\cA]$ whose entries 
and eigenvalues are all in $\rats$: the Bose-Mesner algebra over $\rats$ of the 
subscheme of $\cA$ corresponding to the equivalence classes of $\sim$. 
In this work, we consider matrices in $L[\cA]$ with entries and eigenvalues in 
$\rats(\sqrt{\Delta})$ for some square-free integer $\Delta$. 
More precisely,  
we find all oriented Cayley graphs or signed Cayley graphs $X$ on $G$ 
such that all eigenvalues of $A(X)$ are of the form $m_r\sqrt{\Delta}$
for some fixed square-free integer $\Delta$ and integer $m_r$ 
(with $\Delta < 0$ for oriented Cayley, 
and $\Delta >0$ for signed Cayley).

In this section we find all possible $\Delta$s,  
which are determined by the exponent of the abelian group $G$. 
To do this,  
we find all the quadratic subfields of 
the cyclotomic field $\rats(\zeta_{n})$ for any positive integer $n$. 
The case when $n$ is an odd prime is well-known (see Lemma~\ref{lem:rootpprime}). 
We use the following result from Galois theory.
 
\begin{theorem}\label{thm:Galois}\cite{galoiscox}
Let $E$ be the splitting field of a separable polynomial over a field $F$, 
and let $\Ga=\Gal(E/F)$. 
Then
\begin{enumerate}
\item[(a)]
There is a one-to-one order-reversing correspondence between the subgroups 
of $\Ga$ and the subfields of $E$ that contain $F$:
\begin{itemize}
\item 
If $H$ is a subgroup of $\Ga$, then the corresponding subfield 
is $E^H=\{a\in E\,|\, \theta(a)=a \text{ for all } \theta\in H\}$, 
and $H=\Gal(E/E^H)$.
\item
If $K$ is a subfield of $E$ that contains $F$, 
then the corresponding subgroup of $\Ga$ is $\Gal(E/K)$, and 
$K=E^{\Gal(E/K)}.$
\end{itemize}
\item[(b)] 
For any subgroup $H$ of $\Ga$, 
\[
[E:E^H]=|H| \quad \text{ and } \quad [E^H:F]=[\Ga:H].\qed
\]
\end{enumerate}
\end{theorem} 

This theorem tells us that to count quadratic subfields of $\rats(\zeta_{n})$ over $\rats$, 
we just need to find the number of subgroups of $\Gal(\rats(\zeta_{n})/\rats)$ of index 2. 

In later sections, we will need to write $\sqrt{\Delta}$ as a $\rats$-linear 
combination of $n$-th roots of unity, for which we introduce
quadratic residues and the Legendre symbol. 
Let $p$ be an odd prime. 
An integer $k$ is called a \textsl{quadratic residue} modulo $p$ if $k \equiv x^2 \pmod p$ 
for some $x$; otherwise it is called a \textsl{quadratic non-residue} modulo $p$. 
The \textsl{Legendre symbol} is a function of $k$ and $p$: 
\[
\left(\frac k p\right) = 
\begin{cases*}
      1, & \text{ if $k$ is a quadratic residue modulo $p$, 
      	and  $k\not\equiv 0\pmod p$}; \\
     -1,   & \text{ if $k$ is a quadratic non-residue modulo $p$;}\\
      0, & \text{ if  $k \equiv 0 \pmod p$}.
    \end{cases*}
\]
For a general positive integer $n$, a similar notion for $k\in \Z_n^*$ and $n$  
is defined, with the corresponding symbol called the \textsl{Jacobi symbol}.

Let $p$ be an odd prime. 
The following result of writing $\sqrt{(-1)^{\frac{p-1}{2}}p}$ as a linear combination 
of primitive $p$-th roots of unity plays an important role in finding a similar 
expression for a general square-free integer $\Delta$ in the next section. 

\begin{lemma}\cite{Lemmermeyer}\label{lem:rootpprime}
Let $p$ be an odd prime. 
Then there is a unique quadratic subfield of the cyclotomic field $\rats(\zeta_p)$.
The unique quadratic subfield is $\rats(\sqrt{p})$ if $p \equiv 1\pmod 4$  
and is $\rats(\sqrt{-p})$ if $p \equiv 3\pmod 4$. 
Furthermore, 
\[
\sum_{j=1}^{p-1}\left(\frac jp\right)\zeta_p^j =\sqrt{(-1)^{(p-1)/2}p}=
\begin{cases}
\sqrt{p} & \text{ if $p\equiv 1 \pmod 4$},\\
\sqrt{-p} & \text{ if $p\equiv 3 \pmod 4$}.
\end{cases}
\]
\end{lemma}
We also need the following result on the number of subgroups of index 2 of a given group:
\begin{lemma}\cite{index2subg} \label{lem:subgind2}
Let $n$ be a nonnegative integer and $G$ be the direct product of $n$ cyclic groups each of even order. 
Then the number of distinct subgroups of index 2 in the group $G$ is $2^n-1$.\qed
\end{lemma}
Now we count the number of quadratic subfields of a cyclotomic field. 
\begin{lemma}\label{lem:subind2cyclo}
Let $n=2^mp_1^{m_1}\cdots p_r^{m_r}$ be the prime decomposition of $n>1$. 
Then the number $s$ of quadratic subfields of $\rats(\zeta_{n})$ satisfies
\[
    s=
\begin{cases}
      2^{2+r}-1 & \text{ if } m\geq 3, \\
      2^{1+r}-1  & \text{ if } m=2,\\
      2^r-1 & \text{ if  } m=0,1.
    \end{cases}
\]
\end{lemma}
\proof
By Theorem~\ref{thm:Galois}, 
we just need to count the number subgroups of 
$\Gal(\rats(\zeta_{n})/\rats)$ of index 2. 
It is known that 
\[
 \Z_{2^m}^*\cong 
\begin{cases}
\Z_2 \times \Z_{2^{m-2}} & \text{ if } m\geq 3,\\
\Z_2 & \text{ if }m=2,\\
\Z_1 &\text{ if }m=1.
\end{cases} 
\]
Now the result follows from 
\[
\Gal(\rats(\zeta_{n})/\rats)\cong \Z_{n}^*\cong 
\Z_{2^m}^* \times \Z_{p_1^{m_1}}^*\times \cdots \times \Z_{p_r^{m_r}}^*, 
\]
the fact $\Z_{p^k}^*$ is cyclic of even order for any odd prime $p$ and Lemma~\ref{lem:subgind2}.
\qed

Next we find all square-free integers $\Delta$ such that $\rats(\sqrt{\Delta})$ is a quadratic subfield of $\rats(\zeta_{n})$.

\begin{lemma}\label{lem:deltalist}
Let $n>1$ be a positive integer with prime decomposition $n=2^mp_1^{m_1}\ldots p_r^{m_r}$.  
Let $\Delta\neq1$ be a square-free integer (positive or negative). 
Then $\rats(\sqrt{\Delta})$ is a quadratic subfield of $\rats(\zeta_{n})$ 
if and only if 
\begin{equation}\label{eq:deltalist}
\Delta 
=
2^{w_0} 
(-1)^w
\big((-1)^{(p_1-1)/2}p_1\big)^{\! w_1}\cdots 
\big((-1)^{(p_r-1)/2}p_r\big)^{\! w_r},
\end{equation}
for some $w,w_0,w_1,\ldots w_r\in \{0,1\}$  
such that  $w_0=0$ if $m<3$ and $w=0$ if $m<2$. 
\end{lemma}

\proof 
Let $L=\rats(\zeta_{n})$. 
If $n=2$, 
then $L=\rats$, 
and it contains no quadratic subfield.
If $n=4$, 
then $L=\rats(i)$. 
The only quadratic subfield of $L$ over $\rats$ is $L$ itself, 
and 
\[
2i=i^1-i^3=\left(\frac {1}{4}\right)\zeta_4+\left(\frac {3}{4}\right)\zeta_4^3, 
\]
where $\left(\frac {a}{4}\right)$ is the Jacobi symbol.\\
If $n=8$, 
then $L=\rats(\zeta_8)$, 
and $\Gal(\rats(\zeta_{8}) / \rats)\cong \Z_{8}^*
\cong \Z_2\times \Z_2$.
The three quadratic subfields of $L$ are $\rats(i)$, 
$\rats(\sqrt{2})$ and $\rats(\sqrt{-2})$. 
Furthermore, 
\[
2\sqrt{2}=\zeta_8^1
-\zeta_8^3
-\zeta_8^5
+\zeta_8^7,
\text{ and } \;
2\sqrt{-2}=\zeta_8^1
+\zeta_8^3
-\zeta_8^5
-\zeta_8^7.
\] 
Now we consider a general positive integer $n$. 
By Lemma~\ref{lem:rootpprime}, 
we know for any odd prime divisor $p$ of $n$, 
\[
\sqrt{(-1)^{(p-1)/2}p}\in \rats(\zeta_{p}) \subseteq \rats(\zeta_{n}). 
\]
Note 
for two distinct square-free integers $\Delta_1\neq \Delta_2$, 
if $\rats(\sqrt{\Delta_1})$ and $\rats(\sqrt{\Delta_2})$
are quadratic subfields of $L$, 
then so is $\rats(\sqrt{\Delta_1\Delta_2})$. 
Hence $\rats(\sqrt{\Delta})$ is a subfield of $\rats(\zeta_n)$ if 
$\Delta$ is obtained by multiplying some of the factor 
 $(-1)^{(p-1)/2}p$ for any odd prime divisor $p$ of $n$, and possibly multiplied by some of $ -1$,  $2$, $-2$ depending on the multiplicity of 2 as a factor of $n$. 
By the counts of quadratic subfields of $\rats(\zeta_n)$ in Lemma~\ref{lem:subind2cyclo}, 
these are all the possible $\Delta$s.\qed

\section{$\sqrt{\Delta}$ as a linear combination of roots of unity}
\label{sec:getdelta}

Let $n$ be a positive integer. 
In the previous section we have found all the 
square-free integers $\Delta$ such that 
$\rats(\sqrt{\Delta})$ is a quadratic subfield of $\rats(\zeta_{n})$. 
In this section, 
we find a formula to write each of such $\sqrt{\Delta}$
as a linear combination of $n$-th roots of unity. 

We will make use of the following identities related to 
$\sqrt{(-1)^{(p-1)/2}p}\in\rats(\zeta_{p^m})$, 
$\sqrt{2}, \sqrt{-2}\in \rats(\zeta_{2^m})$ for $m\geq 3$,
and $\sqrt{-1}\in  \rats(\zeta_{2^m})$ for $m\geq 2$.

\begin{lemma}\label{lem:sqrtprimepower}
Let $p$ be an odd prime, 
$m$ be a positive integer, 
and $w\in\{0,1\}$. 
Let $\chi$ be the character on $\Z_8^*$ such that $\chi(1)=\chi(7)=1, \chi(3)=\chi(5)=-1$. 
Let $\left(\frac {k}{n}\right)$ denote the Legendre or Jacobi symbol.  
Then 
\begin{itemize}
\item
For any integer $m\geq 2$ and any odd integer $t$, 
\begin{equation}\label{eq:2^mimp}
\left(\frac {t}{4}\right)2^{m-1}\sqrt{-1} 
=\sum_{\substack{j=1\\(j,2)=1}}^{2^m-1}
\left(\frac {j}{4}\right)\zeta_{4}^{tj}.
\end{equation}
\item
For any integer $m\geq 3$ and any odd integer $t$, 
\begin{equation}\label{eq:2^msqrt2mp}
\chi(t)\left(\frac {t}{4}\right)^{\!\!w}
2^{m-2}\sqrt{2}(\sqrt{-1})^w
=\sum_{\substack{j=1\\(j,2)=1}}^{2^m-1}
\chi(j)\left(\frac {j}{4}\right)^{\!\!w}
\zeta_{8}^{tj}.
\end{equation}
\item
For any positive integer $m$ and any integer $t$ not divisible by $p$, 
\begin{equation} \label{eq:p^msqrtpmp}
\left(\frac {t}{p}\right)^{\!\!w}(-1)^{w+1}p^{m-1}\left(\sqrt{(-1)^{(p-1)/2}p}\right)^{\!\!w}
=\sum_{\substack{j=1 \\ (j,p)=1}}^{p^m-1}
\left(\frac {j}{p}\right)^{\!\!w}\zeta_{p}^{tj} .
\end{equation}
\end{itemize}
\end{lemma}

\proof
Assume $m\geq 2$ and $(t,2)=1$. 
Then by use of  the multiplicative property of Jacobi symbol and that 
$\left(\frac {t^{-1}}{4}\right)=\left(\frac {t}{4}\right)$ for any $t$ such that $(t,4)=1$,  
we have 
\begin{align*}
\sum_{\substack{j=1\\(j,2)=1}}^{2^m-1}
\left(\frac {j}{4}\right)\zeta_{4}^{tj} 
& =\sum_{\ell \in \Z_4^*} \sum_{k=0}^{2^{m-2}-1}
\left(\frac {4k+\ell}{4}\right)\zeta_{4}^{t(4k+\ell)} \quad (\text{let }j=4k+\ell)\\
& = \left(\sum_{\ell \in \Z_4^*}\left(\frac {\ell}{4}\right)\zeta_{4}^{t\ell}\right)
( \sum_{k=0}^{2^{m-2}-1}1)\\
& = \left(\sum_{\ell' \in \Z_4^*}\left(\frac {t^{-1}\ell'}{4}\right)\zeta_{4}^{\ell'}\right)
(2^{m-2}) \quad (\text{let } t\ell =\ell') \\
& =\left(\frac {t}{4}\right) (\sum_{\ell' \in \Z_4^*}\left(\frac {\ell'}{4}\right)\zeta_{4}^{\ell'})
(2^{m-2}) \\
& = \left(\frac {t}{4}\right)(2i)(2^{m-2})\\
& = \left(\frac {t}{4}\right)2^{m-1}i.
\end{align*}
By use of 
\[
\sum_{j\in \Z_8^*}\chi(j)\zeta_{8}^{j} =2\sqrt{2}, \;
\sum_{j\in \Z_8^*}\chi(j)\left(\frac {j}{4}\right)\zeta_{8}^{j} =2\sqrt{-2}, 
\]
and  
\[
\sum_{j\in \Z_p^*}\left(\frac {j}{p}\right)\zeta_{p}^{j} =\sqrt{(-1)^{(p-1)/2}p}, \;
\sum_{j\in \Z_p^*}\zeta_{p}^{j}=-1,
\]
the other two equations follow similarly.\qed

Now we make use of Chinese Remainder Theorem to write $\sqrt{\Delta}$, for 
each of the $\Delta$s in \eqref{eq:deltalist} (so $\rats(\sqrt{\Delta})$ is a subfield of $\rats(\zeta_n)$), as a linear combination of $n$-th roots of unity with coefficients $\pm a\in\rats$, where $a$ is determined by $n$ and $\Delta$. 

\begin{theorem}[Chinese Remainder Theorem] \label{thm:chineserem}
Let $n_1,\ldots, n_k>1$ be a set of pairwise coprime integers and  
let $N=n_1\cdots n_k$.  
Let $a_1,\ldots, a_k$ be any integers.  
Then the system
\begin{align*}
x & \equiv a_1 \pmod{n_1}\\
& \vdots \\
x & \equiv a_k \pmod{n_k},
\end{align*}
has a solution, 
and any two solutions,
say $x_1$ and $x_2$,
satisfies $x_1\equiv x_2 \pmod{N}$.

Furthermore, 
for each $j=1,\ldots, k$, 
let $N_j=N/n_j$,  
and let $s_j \in \Z_{n_j}$ such that $N_js_j \equiv 1 \pmod{n_j}$. 
Then 
\[
x=\sum_{j=1}^k a_j s_j N_j
\]
is a solution of the system. \qed
\end{theorem}
 
\medskip

\begin{lemma}\label{lem:getsqrtDelta}
Let $n=2^mp_1^{m_1}\cdots p_r^{m_r}$ be the prime decomposition of an integer $n>1$. 
Let $w,w_1,\ldots, w_r\in \{0,1\}$.  
Let  
\allowdisplaybreaks
{\footnotesize
\begin{align*}
&a_1=2^{m-2}(-1)^{w_1+\cdots+w_r+r}p_1^{m_1-1}\cdots p_r^{m_r-1}
\chi(p_1\cdots p_r)
\left(\frac {p_1\cdots p_r}{4}\right)^{\!w}
\left(\frac {2p_2\cdots p_r}{p_1}\right)^{\!\!w_1} \cdots
\left(\frac {2p_1\cdots p_{r-1}}{p_r}\right)^{\!\!w_r},\\
&a_2=2^{m-1}(-1)^{w_1+\cdots+w_r+r}p_1^{m_1-1}\cdots p_r^{m_r-1}
\left(\frac {p_1\cdots p_r}{4}\right)
\left(\frac {p_2\cdots p_r}{p_1}\right)^{\!\!w_1} \cdots
\left(\frac {p_1\cdots p_{r-1}}{p_r}\right)^{\!\!w_r},\\
&a_3=2^{m-1}(-1)^{w_1+\cdots+w_r+r}p_1^{m_1-1}\cdots p_r^{m_r-1}
\left(\frac {p_2\cdots p_r}{p_1}\right)^{\!\!w_1} \cdots
\left(\frac {p_1\cdots p_{r-1}}{p_r}\right)^{\!\!w_r}, \text{ and}\\
&a_4=(-1)^{w_1+\cdots+w_r+r}p_1^{m_1-1}\cdots p_r^{m_r-1}
\left(\frac {p_2\cdots p_r}{p_1}\right)^{\!\!w_1} \cdots
\left(\frac {p_1\cdots p_{r-1}}{p_r}\right)^{\!\!w_r}.\end{align*}
}
If $m\geq 3$, then
\begin{align*}
& a_1
\sqrt{2
(-1)^w
[(-1)^{(p_1-1)/2}p_1]^{w_1}\cdots 
[(-1)^{(p_r-1)/2}p_r]^{w_r}}\\
& =\sum_{\substack{j=1\\(j,n)=1}}^{n-1}
\chi(j)
\left(\frac {j}{4}\right)^{\!\!w}
\left(\frac {j}{p_1}\right)^{\!\!w_1}\cdots
\left(\frac {j}{p_r}\right)^{\!\!w_r}
\zeta_{n}^{\frac{n}{8p_1\cdots p_r}j}.
\end{align*}
If $m\geq 2$, then 
\begin{align*}
& a_2
\sqrt{
(-1)
[(-1)^{(p_1-1)/2}p_1]^{w_1}\cdots 
[(-1)^{(p_r-1)/2}p_r]^{w_r}}\\
& =\sum_{\substack{j=1\\(j,n)=1}}^{n-1}
\left(\frac {j}{4}\right)
\left(\frac {j}{p_1}\right)^{\!\!w_1}\cdots
\left(\frac {j}{p_r}\right)^{\!\!w_r}
\zeta_{n}^{\frac{n}{4p_1\cdots p_r}j}.
\end{align*}
If $m\geq 0$, then for $a=a_3$ if $m\geq 1$ and $a=a_4$ if $m=0$, 
we have 
\begin{align*}
& a
\sqrt{
[(-1)^{(p_1-1)/2}p_1]^{w_1}\cdots 
[(-1)^{(p_r-1)/2}p_r]^{w_r}}\\
& =\sum_{\substack{j=1\\(j,n)=1}}^{n-1}
\left(\frac {j}{p_1}\right)^{\!\!w_1}\cdots
\left(\frac {j}{p_r}\right)^{\!\!w_r}
\zeta_{n}^{\frac{n}{p_1\cdots p_r}j}.
\end{align*}
\end{lemma}

\proof
Assume  $m\geq 3$. 
Recall that 
$\Z_{n}^* \cong 
\Z_{2^m}^* \times \Z_{p_1^{m_1}}^*\times \cdots \times \Z_{p_r^{m_r}}^*$, 
through the isomorphism 
$\psi([j]_{n})=([j]_{2^m},[j]_{p_1^{m_1}},\ldots, [j]_{p_r^{m_r}})$, 
where the subscript of  $[j]_{p_k^{m_k}}$ indicates that the element is in the group 
$\Z_{p_k^{m_k}}^*$.
For $k=1,\ldots, r$, 
let $h_k\in \Z_{p_k^{m_k}}^*$ be the element such that $(n/p_k^{m_k})h_k \equiv 1 \pmod{p_k^{m_k}}$,
and let $h_0\in \Z_{2^{m}}^*$ be the element such that $(n/2^{m})h_0 \equiv 1 \pmod{2^{m}}$.
Then by Chinese Reminder Theorem, for any $j_0\in \Z_{2^m}^*, j_1\in \Z_{p_1^{m_1}}^*, \ldots, j_r\in \Z_{p_r^{m_r}}^*$, the element 
\[
j =j_0h_0(n/2^m)+ j_1h_1(n/p_1^{m_1})+\cdots + j_rh_r(n/p_r^{m_r}) \pmod{n}
\]
is the preimage of 
$([j_0]_{2^m},[j_1]_{p_1^{m_1}},\ldots, [j_r]_{p_r^{m_r}})$ 
under $\psi$ in $\Z_{n}^*$.
Note that 
\[
j\equiv j_0 \pmod{2^m},
\left(\frac {j}{4}\right)=\left(\frac {j_0}{4}\right),  \text{ and } \chi(j)=\chi(j_0),
\]
\[
j\equiv j_k \pmod{p_k^{m_k}} \text{ and }
\left(\frac {j}{p_k}\right)=\left(\frac {j_k}{p_k}\right), \, k = 1,\ldots, r. 
\]
Now we write $\sqrt{\Delta}$ as a linear combination of  $n$-th roots of unity for 
$\Delta=2(-1)^w[(-1)^{(p_1-1)/2}]^{w_1}\cdots [(-1)^{(p_r-1)/2}]^{w_r}$
when $8\mid n$. 

\allowdisplaybreaks
{\small
\begin{align}
& \sum_{\substack{j=1\\(j,n)=1}}^{n-1}
\chi(j)
\left(\frac {j}{4}\right)^{\!\!w}
\left(\frac {j}{p_1}\right)^{\!\!w_1}\cdots
\left(\frac {j}{p_r}\right)^{\!\!w_r}
\zeta_{n}^{\frac{n}{8p_1\cdots p_r}j}\nonumber\\
& = \sum_{\substack{j_0=1 \\ (j_0,2)=1}}^{2^m-1}
        \sum_{\substack{j_1=1 \\ (j_1,p_1)=1}}^{p_1^{m_1}-1} \cdots
         \sum_{\substack{j_r=1 \\ (j_r,p_r)=1}}^{p_r^{m_r}-1}
         \chi(j_0)
         \left(\frac {j_0}{4}\right)^{\!\!w}
  	\left(\frac {j_1}{p_1}\right)^{\!\!w_1}\cdots
	\left(\frac {j_r}{p_r}\right)^{\!\!w_r}
	\zeta_{8p_1\cdots p_r}^{j_0h_0\frac{n}{2^m}+j_1h_1\frac{n}{p_1^{m_1}}+\cdots + j_rh_r\frac{n}{p_r^{m_r}}}\nonumber\\
& = \left(\sum_{\substack{j_0=1 \\ (j_0,2)=1}}^{2^m-1}\chi(j_0)
              \left(\frac {j_0}{4}\right)^{\!\!w}
              \zeta_{8}^{j_0h_0\frac{n}{2^mp_1\cdots p_r}}
        \right)
      \left( \sum_{\substack{j_1=1 \\ (j_1,p_1)=1}}^{p_1^{m_1}-1} 
         \left(\frac {j_1}{p_1}\right)^{\!\!w_1}
         \zeta_{p_1}^{j_1h_1\frac{n}{8p_1^{m_1}p_2\cdots p_r}}
      \right) 
      \cdots   \nonumber\\
 & \;\;\;\;\;    \left( \sum_{\substack{j_r=1 \\ (j_r,p_r)=1}}^{p_r^{m_r}-1} 
         \left(\frac {j_r}{p_r}\right)^{\!\!w_r}
         \zeta_{p_r}^{j_rh_r\frac{n}{8p_1p_2\cdots p_r^{m_r}}}
      \right) . \label{eq:delbreak}
\end{align}
}
Now we evaluate the two types of sums abo ve. 
Since $(h_0\frac{n}{2^mp_1\cdots p_r},2)=1$, 
\begin{align}
&\sum_{\substack{j_0=1 \\ (j_0,2)=1}}^{2^m-1}\chi(j_0)
              \left(\frac {j_0}{4}\right)^{\!\!w}
              \zeta_{8}^{j_0h_0\frac{n}{2^mp_1\cdots p_r}} \nonumber\\
& =\chi(h_0\frac{n}{2^mp_1\cdots p_r})
       \left(\frac {h_0\frac{n}{2^mp_1\cdots p_r}}{4}\right)^{\!\!w}
       2^{m-2}\sqrt{2}(\sqrt{-1})^w \quad \text{ ($t=h_0\frac{n}{2^mp_1\cdots p_r}$ in \eqref{eq:2^msqrt2mp}) }  \nonumber \\
&=   \chi(p_1\cdots p_r) \left(\frac {p_1\cdots p_r}{4}\right)^{\!\!w}     
       2^{m-2}\sqrt{2}(\sqrt{-1})^w,   \label{eq:sqrt2ag} 
\end{align}
where in the last equality we used the multiplicative property of 
$\chi$ and of the Legendre symbol, 
and the fact that $h_0\frac{n}{2^m}\equiv 1\pmod{2^m}$.
Similarly, 
using the fact $(h_k\frac{n}{8p_1\cdots p_{k-1}p_k^{m_k}p_{k+1}\cdots p_r},p_k)=1$ and 
\eqref{eq:p^msqrtpmp} with $t=h_k\frac{n}{8p_1\cdots p_{k-1}p_k^{m_k}p_{k+1}\cdots p_r}$, 
it follows that  
for $k=1,\ldots, r$,
\begin{align}
& \sum_{\substack{j_k=1 \\ (j_k,p_k)=1}}^{p_k^{m_k}-1} 
         \left(\frac {j_k}{p_k}\right)^{\!\!w_k}
         \zeta_{p_k}^{j_kh_k\frac{n}{8p_1\cdots p_{k-1}p_k^{m_k}p_{k+1}\cdots p_r}} \nonumber \\
&=  \left(\frac {2p_1\cdots p_{k-1}p_{k+1}\cdots p_r} {p_k}\right)^{\!\!w_k}  
(-1)^{w_k+1}p_k^{m_k-1}\left(\sqrt{(-1)^{(p_k-1)/2}p_k}\right)^{\!\!w_k} .\label{eq:sqrtpag}
 \end{align}        

Therefore continuing with equation \eqref{eq:delbreak} and 
making use of \eqref{eq:sqrt2ag}  and \eqref{eq:sqrtpag}, 
we have
{\small
\begin{align*}
& \sum_{\substack{j=1\\(j,n)=1}}^{n-1}
\chi(j)
\left(\frac {j}{4}\right)^{\!w}
\left(\frac {j}{p_1}\right)^{\!w_1}\cdots
\left(\frac {j}{p_r}\right)^{\!w_r}
\zeta_{n}^{\frac{n}{8p_1\cdots p_r}j}\\
   & =\chi(p_1\cdots p_r) \left(\frac {p_1\cdots p_r}{4}\right)^{\!w}     
       2^{m-2}\sqrt{2}(\sqrt{-1})^w
        \left(\frac {2p_2\cdots p_r} {p_1}\right)^{\!\! w_1}  
(-1)^{w_1+1}p_1^{m_1-1} \\  
 & \;\;\;\;\;    \left(\sqrt{(-1)^{(p_1-1)/2}p_1}\right)^{\!\! w_1}
  \cdots \left(\frac {2p_1\cdots  p_{r-1}} {p_r}\right)^{\!\!  w_r}  
(-1)^{w_r+1}p_r^{m_r-1}\left(\sqrt{(-1)^{(p_r-1)/2}p_r}\right)^{\!\!  w_r} \\
& = 2^{m-2}(-1)^{w_1+\cdots +w_r+r}p_1^{m_1-1}\cdots p_r^{m_r-1}
       \chi(p_1\cdots p_r) \left(\frac {p_1\cdots p_r}{4}\right)^{ \! w}  
       \left(\frac {2p_2\cdots p_r} {p_1}\right)^{ \!\! w_1} 
       \cdots   \\
& \;\;\;\;\;  \left(\frac {2p_1\cdots  p_{r-1}} {p_r}\right)^{\!\!  w_r}
\sqrt{2(-1)^w\left[(-1)^{(p_1-1)/2}p_1\right]^{w_1}\cdots
		      \left[(-1)^{(p_r-1)/2}p_r\right]^{w_r}},
\end{align*}
}
proving the case when $m\geq 3$. 
The other three cases follows similarly. 
We can put all the cases together as
\begin{align}
& \sum_{\substack{j=1\\(j,n)=1}}^{n-1}
\chi(j)^{w_0}
\left(\frac {j}{4}\right)^{\!\!w}
\left(\frac {j}{p_1}\right)^{\!\!w_1}\cdots
\left(\frac {j}{p_r}\right)^{\!\!w_r}
\zeta_{n}^{\frac{n}{\min\{4^w8^{w_0},8\}p_1\cdots p_r}j} \nonumber \\
& =a
\sqrt{2^{w_0}
(-1)^{w}
[(-1)^{(p_1-1)/2}p_1]^{w_1}\cdots 
[(-1)^{(p_r-1)/2}p_r]^{w_r}}, \label{eq:comsqrtdel} 
\end{align}
where $w_0=0$ if $m<3$,
$w=0$ if $m<2$, and \\
$
a=(-1)^{w_1+\cdots+w_r+r}p_1^{m_1-1}\cdots p_r^{m_r-1} 
2^{\max\{m,1\}-1-w_0}
\chi(p_1\cdots p_r)^{w_0}
\left(\frac {p_1\cdots p_r}{4}\right)^{\!w} $\\
$
\left(\frac {2^{w_0}p_2\cdots p_r}{p_1}\right)^{\!w_1} \cdots
\left(\frac {2^{w_0}p_1\cdots p_{r-1}}{p_r}\right)^{\!w_r}.
$
\qed

\section{Characterization of the graphs}
\label{sec:mainthm}
In this section, 
we consider oriented and signed Cayley graphs with desired eigenvalue form. 
Here, 
when talking about adjacency matrix of an oriented graph $X$ on $n$ vertices, 
we mean
the $(1,-1,0)$ skew-symmetric matrix of size $n\times n$ defined by:
 \[
A(X)_{u,v} = 
\begin{cases*}
      1 & \text{ if $(u,v)\in \arc(X)$}, \\
     -1   & \text{ if $(v,u)\in \arc(X)$,}\\
      0 & \text{ otherwise}.
    \end{cases*}
\]
When talking about the adjacency matrix of a signed graph $X$, 
we mean the $(1,-1,0)$ symmetric matrix such that $A_{u,v}=A_{v,u}=1$ if $\{u,v\}$ is an edge of $X$ with a positive sign, 
$A_{u,v}=A_{v,u}=-1$ if $\{u,v\}$ is an edge of $X$ with a negative sign, 
and $A_{u,v}=A_{v,u}=0$ otherwise.

For a finite abelian group $G$,
Corollary~\ref{thm:integtran} tells us that
there are Cayley graphs $X(G,C)$ all of whose eigenvalues are integers,
and a characterization of $C$ for $X(G,C)$ being such a graph is given. 
Denote the group scheme of $G$ by $\cA$ and assume the splitting field of $\cA$ is $L$. 
Here we consider when a $(0,1,-1)$ matrix in the Bose-Mesner algebra $L[\cA]$ has all its eigenvalues 
integer multiples of $\sqrt{\Delta}$ for some square-free 
integer $\Delta\neq1$.
When $\Delta>1$, 
we have a signed Cayley graph;
when $\Delta<0$, 
we have an oriented Cayley graph. 
For all the $\Delta$s in Lemma~\ref{lem:deltalist}, 
there exists a corresponding signed or oriented Cayley graph. 
We give a characterization of the connection set of such Cayley graphs. 
 
For a subset $C$ of $G$, 
we use $A_C$ to denote the matrix $\sum_{g\in C}A_g$.  
Let $n_1$ be a positive integer. 
Recall that $\Gal(\rats(\zeta_{n_1})/\rats)\cong \Z_{n_1}^*$, 
with each $k\in \Z_{n_1}^*$ corresponding to $\tau_k$, 
where $\tau_k(\zeta_{n_1})=\zeta_{n_1}^k$. 

\begin{theorem} \label{thm:char}
Let $G$ be a finite abelian group with exponent $n_1$. 
Assume $n_1=2^mp_1^{m_1}\cdots p_r^{m_r}$ 
is the prime decomposition of $n_1$. 
Let $\chi$ be the character on $\Z_8^*$ such that 
$\chi(1)=\chi(7)=1$ and $\chi(3)=\chi(5)=-1$. 
Let 
\begin{equation}\label{eq:deltas}
\Delta 
=
2^{w_0} 
(-1)^w
((-1)^{(p_1-1)/2}p_1)^{w_1}\cdots 
((-1)^{(p_r-1)/2}p_r)^{w_r}
\end{equation}
for some $w, w_0, w_1,\ldots, w_r\in\{0,1\}$ such that 
$w_0=0$ if $m<3$ and $w=0$ if $m<2$. 
Let $H_{\Delta}$ be the subgroup of $\Z_{n_1}^*$ defined by 
\begin{equation} \label{eq:HFsqrtDel}
H_{\Delta}=\{k\in \Z_{n_1}^*\,|\, 
\chi(k)^{w_0}
\left(\frac {k}{4}\right)^{\!\!w}
\left(\frac {k}{p_1}\right)^{\!\!w_1}\cdots
\left(\frac {k}{p_r}\right)^{\!\!w_r}=1
\}. 
\end{equation}
For each of the above $\Delta$, define an equivalence relation $\sim_\Delta$ on $G$ by $g\sim_{\Delta} h$ 
if and only if $g=kh$ for some $k\in H_{\Delta}$. \\
Denote the group scheme of $G$ by $\cA$. 
Let $\Delta_0$ be a square-free integer.  
Then there exists a (0,1,-1) matrix in $\cx[\cA]$ all of whose eigenvalues are in 
$\rats(\sqrt{\Delta_0})$ if and only if $\Delta_0$ satisfies \eqref{eq:deltas} for some $w_j$s. Take any such $\Delta$. 
Assume $M=A_{C_1}-A_{C_2}$ for two disjoint subsets $C_1$ and $C_2$ of $G$ such that $0\notin C_1\cup C_2$. 
Then 
\begin{enumerate}
\item[(i)]
All eigenvalues of $M$ lie in $\rats(\sqrt{\Delta})$ if and only if  both 
$C_1$ and $C_2$ are unions of equivalence classes of $\sim_{\Delta}$. Furthermore, 
\item[(ii)]
All eigenvalues of $M$ are integer multiples of $\sqrt{\Delta}$ if and only if either 
\begin{itemize}
\item
$\Delta=1$, and both 
$C_1$ and $C_2$ are unions of equivalence classes of $\sim_{1}$, or 
\item 
$\Delta\neq 1$,  and $C_2=kC_1$ for some $k \in  \Z_{n_1}^* \backslash H_\Delta$ and $C_1$ (hence $C_2$ as well) is a union of equivalence classes of $\sim_\Delta$.  
\end{itemize}
Note that $\Delta<0$ corresponds to $M$ being the adjacency matrix of the oriented Cayley $X(G,C_1)$, and $\Delta >0$ corresponds to $M$ being the adjacency matrix of the signed Cayley $X(G,C_1, C_2)$. 
\end{enumerate}
\end{theorem}
\proof
Since the exponent of $G$ is $n_1$, 
by the fundamental theorem of abelian groups, 
$G \cong \Z_{n_1}\times \Z_{n_2}\times \cdots \times \Z_{n_e}$ 
for some positive integers $n_1,\ldots, n_e$ 
with $n_j\mid n_{j-1}$ for $j=2,\ldots, e$. 
The splitting field of the group scheme $\cA$ over $G$
is $L=\rats(\zeta_{n_1})$, 
and the Krein field is $K=\rats$.
As we have shown in Lemma~\ref{lem:deltalist},
any quadratic subfield $\rats(\sqrt{\Delta})$ 
of  $\rats(\zeta_{n_1})$  satisfies
\[
\Delta 
=
2^{w_0} 
(-1)^w
[(-1)^{(p_1-1)/2}p_1]^{w_1}\cdots 
[(-1)^{(p_r-1)/2}p_r]^{w_r}
\]
for some $w,w_0,w_1,\ldots, w_r\in \{0,1\}$, 
with $w_0=0$ if $m<3$ and $w=0$ if $m<2$.

Take any such $\Delta$. 
Let $F=\rats(\sqrt{\Delta})$, 
and let $\sim_F$ and $Z_F$ be defined as in Theorem~\ref{thm:HFconnectionset}.  
By Theorem~\ref{thm:HFconnectionset}, 
all the eigenvalues of $M=A_{C_1}-A_{C_2}$ belong to $\rats(\sqrt{\Delta})$ if and only if 
both $C_1$ and $C_2$ are unions of equivalence classes of $\sim_F$.  
To finish proving (1), 
we just need to prove $Z_F=H_{\Delta}$.   
Recall that $Z_F=\{k\in \Z_{n_1}^*\;|\; \tau_k\in \Gal\big(\rats(\zeta_{n_1})/\rats(\sqrt{\Delta})\big)\}$, 
hence we need to prove 
$\{k\in \Z_{n_1}^* \, | \, \tau_k(\sqrt{\Delta})=\sqrt{\Delta}\} =H_\Delta$. 
From equation \eqref{eq:comsqrtdel},
we know that 
\[
\sum_{\substack{j=1\\(j,n_1)=1}}^{n_1-1}
\chi(j)^{w_0}
\left(\frac {j}{4}\right)^{\!\!w}
\left(\frac {j}{p_1}\right)^{\!\!w_1}\cdots
\left(\frac {j}{p_r}\right)^{\!\!w_r}
\zeta_{n_1}^{\frac{n_1}{\min\{4^w 8^{w_0},8\}p_1\cdots p_r}j}
=a
\sqrt{\Delta}
\]
for some nonzero integer $a$. 
Now let $k\in \Z_{n_1}^*$, 
then 

\begin{align*}
&\tau_k(a
\sqrt{\Delta})\\
= &\tau_k(\sum_{\substack{j=1\\(j,n_1)=1}}^{n_1-1}
\chi(j)^{w_0}
\left(\frac {j}{4}\right)^{\!\!w}
\left(\frac {j}{p_1}\right)^{\!\!w_1}\cdots
\left(\frac {j}{p_r}\right)^{\!\!w_r}
\zeta_{n_1}^{\frac{n_1}{\min\{4^w 8^{w_0},8\}p_1\cdots p_r}j})\\
= &\sum_{\substack{j=1\\(j,n_1)=1}}^{n_1-1}
\chi(j)^{w_0}
\left(\frac {j}{4}\right)^{\!\!w}
\left(\frac {j}{p_1}\right)^{\!\!w_1}\cdots
\left(\frac {j}{p_r}\right)^{\!\!w_r}
\zeta_{n_1}^{\frac{n_1}{\min\{4^w 8^{w_0},8\}p_1\cdots p_r}jk}\\
= & \sum_{\substack{j'=1\\(j',n_1)=1}}^{n_1-1}
\chi(j'k^{-1})^{w_0}
\left(\frac {j'k^{-1}}{4}\right)^{\!\!w}
\left(\frac {j'k^{-1}}{p_1}\right)^{\!\!w_1}\cdots
\left(\frac {j'k^{-1}}{p_r}\right)^{\!\!w_r}
\zeta_{n_1}^{\frac{n_1}{\min\{4^w 8^{w_0},8\}p_1\cdots p_r}j'} \text{ (let }j'=jk)\\
 = & \chi(k)^{w_0}
\left(\frac {k}{4}\right)^{\!\!w}
\left(\frac {k}{p_1}\right)^{\!\!w_1}\cdots
\left(\frac {k}{p_r}\right)^{\!\!w_r}
\sum_{\substack{j=1\\(j,n_1)=1}}^{n_1-1}
\chi(j)^{w_0}
\left(\frac {j}{4}\right)^{\!\!w} \\
& \;\;
\left(\frac {j}{p_1}\right)^{\!\!w_1}\cdots
\left(\frac {j}{p_r}\right)^{\!\!w_r}
\zeta_{n_1}^{\frac{n_1}{\min\{4^w 8^{w_0},8\}p_1\cdots p_r}j}\\
= &\chi(k)^{w_0}
\left(\frac {k}{4}\right)^{\!\!w}
\left(\frac {k}{p_1}\right)^{\!\!w_1}\cdots
\left(\frac {k}{p_r}\right)^{\!\!w_r} 
a\sqrt{\Delta}. 
\end{align*}
Therefore $\tau_k\in \Gal(\rats(\zeta_{n_1})/\rats)$ fixes 
$\sqrt{\Delta}$, 
that is, $\tau_k\in \Gal(\rats(\zeta_{n_1})/\rats(\sqrt{\Delta}))$, 
if and only if 
\[
\chi(k)^{w_0}
\left(\frac {k}{4}\right)^{w}
\left(\frac {k}{p_1}\right)^{w_1}\cdots
\left(\frac {k}{p_r}\right)^{w_r}=1.
\]
We have proved (i). 

Now we prove (ii). 
Since $M=A_{C_1}-A_{C_2}$ is an integer matrix, 
all its eigenvalues are algebraic integers. 
Therefore its rational eigenvalues are all integers. 
The case corresponds to $\Delta=1$ follows from (i). \\ 
Assume $\Delta\neq 1$. 
Then $H_\Delta$ is a subgroup of $\Z_{n_1}^*=H_1$ of index 2. 
We first prove the condition is sufficient. 
Assume $C_2=k_0C_1$ for some $k_0\in \Z_{n_1}^*\backslash H_\Delta$. 
Then in fact, 
$C_2=kC_1$ for any $k\in \Z_{n_1}^*\backslash H_\Delta$, 
as $\Z_{n_1}^*$ is the disjoint union of  
$H_\Delta$ and  $k_0H_\Delta$. 
Further assume $C_1$ is a union of $\sim_\Delta$ equivalence classes. 
Then so is $C_2$. 
By part (i), 
all the eigenvalues of $A_{C_j}$ are in $\rats(\sqrt{\Delta})$ for $j=1,2$. 
Since $A_{C_1}$ is an integer matrix, 
all its eigenvalues are algebraic integers and are of the form 
$\frac{a+b\sqrt{\Delta}}{2}$ for some integers $a$ and $b$. 
Now we prove that the eigenvalues of $A_{C_1}$ and of $A_{C_2}$ are related to each other:
on each principal idempotent $E$ of the scheme, 
the corresponding two eigenvalues are algebraic conjugates of each other.  
Label the principal idempotents of the group scheme $\cA$ of $G$ by $g\in G$. 
Assume 
\begin{equation}\label{eq:ac1}
A_{C_1}=\sum_{g\in G}\frac{a_g+b_g\sqrt{\Delta}}{2}E_g.
\end{equation}
Take any $k\in \Z_{n_1}^*\backslash H_\Delta$, 
then $\left(\sqrt{\Delta}\right)^{\tau_k}=-\sqrt{\Delta}$,
and therefore (the following actions of $\tau_k$ and $\htau_k$ on $L[\cA]$ are defined in Section~\ref{mungal}) 
\begin{align*}
& A_{C_1}^{\htau_k}=\sum_{g\in G}\frac{a_g+b_g\sqrt{\Delta}}{2}E_g^{\tau_k}
=\sum_{g\in G}\left(\frac{a_g-b_g\sqrt{\Delta}}{2}\right)^{\!\!\tau_k}E_g^{\tau_k} \\
&=\left(\sum_{g\in G}\frac{a_g-b_g\sqrt{\Delta}}{2}E_g\right)^{\!\!\tau_k}. 
\end{align*}
On the other hand, 
we have 
 $A_{C_1}^{\htau_k}=A_{k^{-1}C_1}=A_{C_2}$, 
 as $k^{-1}\in  \Z_{n_1}^*\backslash H_\Delta$. 
Combining the two equations, 
we have 
$
A_{C_2}=\left(\sum_{g\in G}\frac{a_g-b_g\sqrt{\Delta}}{2}E_g\right)^{\tau_k}. 
$ 
Applying $\tau_k^{-1}$ to both sides of the equation and making use of the fact that $A_{C_2}$ is a $(0,1)$-matrix,  we have 
\begin{equation}\label{eq:ac2}
A_{C_2}=A_{C_2}^{\tau_k^{-1}}=\sum_{g\in G}\frac{a_g-b_g\sqrt{\Delta}}{2}E_g. 
\end{equation}
It follows from equations \eqref{eq:ac1} and \eqref{eq:ac2} that 
\[
M= A_{C_1}-A_{C_2}=\sum_{g\in G}b_g\sqrt{\Delta}E_g, 
\]
all of whose eigenvalues are integer multiples of $\sqrt{\Delta}$ 
and we have proved the sufficiency of the condition. 

Conversely, 
assume the (0,1,-1) matrix $M=A_{C_1}-A_{C_2}$ has all its eigenvalues of the form $m_r\sqrt{\Delta}$.  
By part (i), 
both $C_1$ and $C_2$ are unions of equivalence classes of $\sim_\Delta$, 
therefore by part (i) again, all the eigenvalues of $A_{C_1}$ and $A_{C_2}$ 
are in $\rats(\sqrt{\Delta})$. 
We just need to prove that 
$C_1=kC_2$
for some (therefore any) $k\in \Z_{n_1}^*\backslash H_\Delta$. 
As in proving the sufficiency of the condition, 
we may assume 
$A_{C_1}$ is given as in \eqref{eq:ac1}. 
Since all the eigenvalues of $M=A_{C_1}-A_{C_2}$ are integer multiple of $\sqrt{\Delta}$,  
we know 
\[
A_{C_2}=\sum_{g\in G}\frac{a_g+c_g\sqrt{\Delta}}{2}E_g 
\]
for some integers $c_g$ with $b_g\equiv c_g \pmod2$. 
Take any  $k\in \Z_{n_1}^*\backslash H_\Delta$, 
and as earlier, we have 
\[ 
A_{k^{-1}C_1}=A_{C_1}^{\htau_k}
=\sum_{g\in G}\frac{a_g-b_g\sqrt{\Delta}}{2}E_g \text{ and }
A_{k^{-1}C_2}=A_{C_2}^{\htau_k}
=\sum_{g\in G}\frac{a_g-c_g\sqrt{\Delta}}{2}E_g.
\]
Therefore 
\[ A_{C_1}+ A_{k^{-1}C_1}=
\sum_{g\in G}a_gE_g
=A_{C_2}+A_{k^{-1}C_2}
\]
has only integer eigenvalues and 
\begin{equation}\label{eq:twocon}
C_1\cup k^{-1}C_1=C_2\cup k^{-1}C_2
\end{equation}
 is a union of equivalence classes of $\sim_1$. 
Let $C$ be an equivalence class of $\sim_1$, 
then $k'C=C$ for any $k'\in \Z_{n_1}^*$. 
Now we prove neither $C_1$ nor $C_2$ contains $C$ as a subset. 
Assume on the contrary and without loss of generality that $C\subseteq C_1$, 
then by 
\eqref{eq:twocon},  
we also have $C\subseteq C_2\cup k^{-1}C_2$. 
Since $C_1\cap C_2=\emptyset$, 
we have $C\subseteq k^{-1}C_2$. 
But this implies 
$ C=kC\subseteq C_2$, 
contradicting $C_1\cap C_2=\emptyset$. 
Therefore for each $\sim_1$ equivalence class $C$, 
$C_1$ (and $C_2$) contains at most one of the two $\sim_\Delta$ classes 
$H_\Delta c=\{kc\, | \, k\in H_\Delta\}$ and $C\backslash (H_\Delta c)$ split from it, 
where $c$ is an element of $C$. 
Note that $C\backslash (H_\Delta c)=k H_\Delta c$ for any $k\in \Z_{n_1}^*\backslash H_\Delta$. 
Take any equivalence class $C_0$ of $\sim_\Delta$ such that 
$C_0\subseteq C_1$, 
then equation \eqref{eq:twocon} and the assumption $C_1\cap C_2=\emptyset$ imply  $C_0\subseteq k^{-1} C_2$ and therefore $kC_1\subseteq C_2$. 
Switching the roles of $C_1$ and $C_2$ and repeat the argument for 
$k^{-1}\in   \Z_{n_1}^*\backslash H_\Delta$,  
 we have $k^{-1} C_2\subseteq C_1$. 
Therefore $kC_1=C_2$, 
and the condition is necessary. 
In fact, 
if $A_{C_1}=\sum_{g\in G}\frac{a_g+b_g\sqrt{\Delta}}{2}E_g$ and  
$C_2=kC_1$ for some $k\in   \Z_{n_1}^*\backslash H_\Delta$, 
then $A_{C_2}=\sum_{g\in G}\frac{a_g-b_g\sqrt{\Delta}}{2}E_g$, 
and $M=A_{C_1}-A_{C_2}=\sum_{g\in G}b_g\sqrt{\Delta}E_g$. 

Finally we prove our claim that $M$ is the adjacency matrix of an oriented Cayley or a signed Cayley, 
depending on $\Delta<0$ or $\Delta>0$. 
Let 
$
\Delta =
2^{w_0} 
(-1)^w
((-1)^{(p_1-1)/2}p_1)^{w_1}\cdots 
((-1)^{(p_r-1)/2}p_r)^{w_r}
$ 
as in \eqref{eq:deltas}. 
Then $\Delta<0$ if and only if 
\[
(-1)^{w+\frac{(p_1-1)w_1+\cdots + (p_r-1)w_r}{2}}=-1.
\]
Using the fact $\left(\frac {-1}{4}\right)=-1$ and that 
for an odd prime $p$, 
$\left(\frac {-1}{p}\right)=1$ if and only if $p\equiv 1\pmod 4$,
we know 
\[
(-1)^{w+\frac{(p_1-1)w_1+\cdots + (p_r-1)w_r}{2}}
=\left(\frac {-1}{4}\right)^{\!\! w}
\left(\frac {-1}{p_1}\right)^{\!\! w_1} \cdots
\left(\frac {-1}{p_r}\right)^{\!\! w_r}.
\]
Hence $\Delta<0$ if and only if 
$\left(\frac {-1}{4}\right)^w
\left(\frac {-1}{p_1}\right)^{\! w_1} \cdots
\left(\frac {-1}{p_r}\right)^{\! w_r}=-1$, 
that is, 
if and only if 
$-1 \in  \Z_{n_1}^*\backslash  H_\Delta$, 
as $\chi(-1)=1$. 
Therefore, 
if $\Delta<0$, 
then $C_2=(-1)C_1$. 
Since $C_1\cap C_2=\emptyset$, 
we know $M=A_{C_1}-A_{C_2}= A_{C_1}-A_{C_1}^T$ is the skew-symmetric adjacency matrix of the oriented Cayley $X(G,C)$. 
If $\Delta>0$, 
then $-1\in H_\Delta$ and $-C_i=C_i$ as $C_i$ is a union of $\sim_\Delta$ equivalence classes. 
Hence $M=A_{C_1}-A_{C_2}$ is the symmetric adjacency matrix of the signed Cayley $X(G,C_1,C_2)$.\qed

Note that for a signed Cayley graph $X(G, C(1),C(-1))$, 
we allow $C(1)$ or $C(-1)$ to empty, 
so a Cayley graph is a special signed Cayley graph with $C(-1)=\emptyset$. 

\begin{example}
Let $G=\Z_8$, 
therefore $n_1=8$
and $\Delta=1,2,-1,-2$. 
We have $H_1=\Z_8^*=\{1,3,5,7\}$. 
From \eqref{eq:HFsqrtDel} we have
\[
H_2=\{1,7\}, \; H_{-1}=\{1,5\},\; H_{-2}=\{1,3\}.
\] 
Therefore all eigenvalues of the Cayley graph $X(\Z_8,\{1,3,5,7\})$ 
are integers; all eigenvalues of the oriented Cayley graphs $X(\Z_8,\{1,5\})$ 
and $X(\Z_8,\{1,2,5\})$ are integer multiples of $i=\sqrt{-1}$; 
all eigenvalues of the oriented Cayley $X(\Z_8,\{1,3\})$ are integer 
multiples of $\sqrt{-2}$; 
all eigenvalues of the signed Cayley $X(\Z_8,\{1,7\},\{3,5\})$ 
are integer multiples of $\sqrt{2}$. 
\end{example}

\section{Nonabelian groups}
\label{sec:nonab}
Here we show that some of our results on Cayley graphs for abelian groups can be extended to normal Cayley graphs for nonabelian groups. 

In this section, $G$ is a group of size $n$, and the product (the group operation) of two elements $g$ and $h$ is denoted by $gh$. The group identity is 1 and the inverse of $g$ is denoted by $g^{-1}$. Two elements $g$ and $h$ are \textsl{conjugate} to each other if 
\[
g = x^{-1}hx \text{ for some } x\in G.
\]
This is an equivalence relation on $G$ and it partitions $G$ into conjugacy classes, which we denote by $K_0,K_1,\ldots, K_d$. 
A Cayley graph $X(G,C)$ is called a \textsl{normal Cayley graph} if the connection set $C$ is a union of conjugacy classes of $G$. 

Under the regular representation of $G$, each group element $g$ corresponds to an $n \times n$ permutation matrix $A_g$  with $(A_g)_{ab}=1$ if and only if $ba^{-1}=g$ (we used this for abelian groups in Section~\ref{sec:BMena}). 
For $i=0,\ldots, d$, define $A_i=\sum_{g\in K_i} A_g$. 
Note that if $G$ is not abelian, then $\{A_g\, | \, g\in G\}$ does not form a (commutative) association scheme. But $\cA=\{A_0,\ldots, A_d\}$ does and $\cx[\cA]$ is isomorphic to the center of the group algebra of $G$. $\cA$ is called the \textsl{conjugacy class scheme} of $G$. 
Denote the irreducible characters of $G$ by $\chi_0=1,\ldots, \chi_d$ (the number of irreducible characters of $G$ is equal to the number of its conjugacy classes), with $\chi_0$ being the principal character. 
Then the principal idempotents $\{E_j\, |\, j=0,\ldots,d\}$ of $\cA$ can be described with these characters as   
\begin{equation}\label{eq:idemnon}
(E_j)_{g,h}=\frac{\chi_j(1)}{n} \chi_j(gh^{-1}). 
\end{equation}
The eigenvalue $p_i(j)$ of $A_i$ on $E_j$ is $\frac{1}{\chi_j(1)}\sum_{g\in K_i}\chi_j(g) =\frac{|K_i|}{\chi_j(1)}\chi_j(g_i)$, where $g_i$ is a representative of $K_i$.

Again let $n_1$ be the exponent of $G$, then by a theorem of Brauer, any the irreducible characters $\chi_j$ of $G$ takes value in $\cL=\rats(\zeta_{n_1})$ \cite{Isaacs}. Further, since the product of the two characters of $G$ is again a character and any character of $G$ is a nonnegative integer linear combination of irreducible characters of $G$, 
we know that all the Krein parameters of $\cA$ lie in $\rats$. 
Now we present a result similar to Theorem~\ref{thm:HFconnectionset}.
 
\medskip

\begin{theorem} \label{thm:HFconnenonabel}  
Let $G$ be a group with exponent $n_1$. Denote its conjugacy scheme by $\cA$. 
Let $\cL=\rats(\zeta_{n_1})$ 
and let $F$, $\phi$, $Z_F$ and $\cF$ be defined as in Theorem~\ref{thm:HFconnectionset}. 
Define an equivalence relation $\sim_F$ on $G$ by $g\sim_F h$ if and only if $g$ and $h^k$ are conjugates for some $k\in Z_F$. 
Then $\cF=F[\cB]$ for a subscheme 
 $\cB$  of $\cA$, 
with each element $B_i$ corresponding to an equivalence class of $G$ under $\sim_F$: $B_{[g]_F}=\sum_{h:h\sim_F g}A_h$. 
 In particular,  a normal Cayley digraph $X(G,C)$  has all its eigenvalues lying in $F$ if and only if the connection set $C$ is a union of equivalence classes of $\sim_F$. 
\end{theorem}
\proof
Since the Krein field is the rational field here, Lemma~\ref{lem:hatpremulti} and Theorem~\ref{eigaut} imply that 
the action of $\Ga=\Gal(\cL/\rats)$ over $\cL[\cA]$, as defined in equation~\eqref{eq:fieldlinact}, preserves matrix multiplication and Schur multiplication. By Lemma~\ref{lem:qlinproshurpro}, for any $\sg\in \Ga$, $\hat{\sg}$ 
permutes the elements of $\cA$ and permutes the elements of $\{E_0,\ldots, E_d\}$, 
and hence for any rational matrix $A\in \cL[\cA]$, $A^{\hat{\sg}}$ is again rational and $(A^{\hat{\sg}})^{\sg^{-1}}=A^{\hat{\sg}}$. 
Now we provide an equivalent way to describe the action of $\hat{\sg}$ on $\cL[\cA]$ in this case. 
Assume $A=\sum_r \theta_r E_r$, then $A^{\hat{\sg}}=\sum_r \theta_r E_r^\sg$ and 
\begin{equation}\label{eq:equhatact}
A^{\hat{\sg}} =(A^{\hat{\sg}})^{\sg^{-1}} =(\sum_r \theta_r E_r^\sg)^{\sg^{-1}}=\sum_r \theta_r ^{\sg^{-1}} E_r. 
\end{equation}
Let $\chi$ be a character of $G$ supported by the representation $\mathfrak{X}$ of dimension $m$. 
Take any $g\in G$ and denote the eigenvalues of $\mathfrak{X}(g)$ by $\lambda_j$ for $j=1,\ldots, m$, then $\lambda_j$ is a power of $\zeta_{n_1}$ and for any $k\in \Z_{n_1}^*$, 
\begin{align}
\chi(g)^{\tau_k}
&= \left(\tr(\mathfrak{X}(g))\right)^{\tau_k}=\big(\sum_{j=1}^m \lambda_j\big)^{\tau_k} = \sum_{j=1}^m \lambda_j^{\tau_k} 
= \sum_{j=1}^m \lambda_j^{k} 
 =\tr(\left(\mathfrak{X}(g)\right)^k) \nonumber \\
&=\tr(\mathfrak{X}(g^k)) =\chi(g^k). \label{eq:chaongroup}
\end{align} 
Let $C$ be a union of conjugacy class of $G$, now we find the image of $A_C$ under $\htau_k$  for $k\in \Z_{n_1}^*$.   
Let $E_r$ be as in equation~\eqref{eq:idemnon},  
then $A_C=\sum_r \frac{1}{\chi_r(1)}\big(\sum_{g\in C}\chi_r(g)\big) E_r$. 
By equations~\eqref{eq:equhatact} and \eqref{eq:chaongroup}, 
\begin{align}
(A_C)^{\htau_k}
& =\sum_r \big(\frac{1}{\chi_r(1)}\sum_{g\in C} \chi_r(g)\big)^{\tau_k^{-1}} E_r
=\sum_r \frac{1}{\chi_r(1)}\big(\sum_{g\in C}\chi_r(g)^{\tau_k^{-1}}\big)  E_r  \nonumber \\
& =\sum_r \frac{1}{\chi_r(1)} \sum_{g\in C} \chi_r(g^{k^{-1}}) E_r =A_{C^{k^{-1}}},  \label{eq:hatset} 
\end{align}
where  $k^{-1}$ is taken in $\Z_{n_1}^*$ and $C^{k^{-1}}=\{g^{k^{-1}}\, | \, g\in C\}$. 
Therefore $A_C$ is fixed by $\htau_k$ if and only if $C=C^{k^{-1}}$, and $A_C$ is fixed by $\{\htau_k\, | \, \tau_k\in \Gal(\cL/F)\}$ if and only if $C$ is a union of $\sim_F$ equivalence classes. Now the result follows from Theorem~\ref{fldsub}. \qed

Now we present some results on normal Cayley graphs with nice eigenvalues, which are similar to Theorem~\ref{thm:char}. 
\begin{theorem} \label{thm:charnonab}
Let $G$ be a finite group with exponent $n_1=2^mp_1^{m_1}\cdots p_r^{m_r}$. 
Let $\chi$, $\Delta$ and $H_\Delta$ be as in Theorem~\ref{thm:char}. 
For each $\Delta$, define an equivalence relation $\sim_\Delta$ on $G$ by $g\sim_{\Delta} h$ 
if and only if $g$ and $h^k$ are conjugate for some $k\in H_{\Delta}$. 
Denote the conjugacy class group scheme of $G$ by $\cA$. 
Let $\Delta_0$ be a square-free integer.  
If there exists a (0,1,-1) matrix in $\cx[\cA]$ all of whose eigenvalues are in 
$\rats(\sqrt{\Delta_0})$,  then $\Delta_0$ satisfies \eqref{eq:deltas} for some $w_j$s. Take any such $\Delta$. 
Assume $M=A_{C_1}-A_{C_2}\in \cx[\cA]$ (so $C_1$ and $C_2$ are closed under conjugate) for two disjoint subsets $C_1$ and $C_2$ of $G$ with $0\notin C_1\cup C_2$. 
Then 
\begin{enumerate}
\item[(i)]
All eigenvalues of $M$ lie in $\rats(\sqrt{\Delta})$ if and only if  both 
$C_1$ and $C_2$ are unions of equivalence classes of $\sim_{\Delta}$. Furthermore, 
\item[(ii)]
All eigenvalues of $M$ are integer multiples of $\sqrt{\Delta}$ if and only if either 
\begin{itemize}
\item
$\Delta=1$, and both 
$C_1$ and $C_2$ are unions of equivalence classes of $\sim_{1}$, or 
\item 
$\Delta\neq 1$,  and $C_2=C_1^k$ for some $k \in  \Z_{n_1}^* \backslash H_\Delta$ and $C_1$ is a union of equivalence classes of $\sim_\Delta$.  
\end{itemize}
Furthermore $\Delta<0$ corresponds to an oriented normal Cayley graph and $\Delta>0$ corresponds to a signed normal Cayley graph. 
\end{enumerate}
\end{theorem}
\proof 
Part (i) follows from Theorem~\ref{thm:HFconnenonabel} and a similar argument as in Theorem~\ref{thm:char}. Now we prove part (ii). 

Assume $\Delta\neq 1$ and $C_1$ and $C_2$ are unions of conjugacy classes of $G$. 
Denote the irreducible characters of $G$ by $\chi_0,\ldots, \chi_d$, 
and the corresponding idempotents of the conjugacy class scheme by $E_0,\ldots, E_d$.  
Then the rest of the proof is similar to that of Theorem~\ref{thm:char}; we just need to replace the $E_g$'s by $E_r$'s, and make use of equation~\eqref{eq:hatset} directly for the image of $A_{C_1}$ under $\htau$, where $C_1$ is a union of conjugacy classes.  
For example, for the necessary condition, assume $C_1=C_2^k$ for some $k\in \Z_{n_1}^*\backslash H_\Delta$  and that $C_1$ is a union of $\sim_\Delta$ equivalence classes. 
By part (i) we can assume without loss of generality that $A_{C_1}=\sum_{r=0}^d \frac{a_r+b_r\sqrt{\Delta}}{2}E_r$. 
By equations~\eqref{eq:hatset} and \eqref{eq:equhatact}, we have $A_{C_2}=A_{C_1^{k^{-1}}}
=(A_{C_1})^{\htau_k}=
\sum_{r=0}^d \frac{a_r-b_r\sqrt{\Delta}}{2}E_r$. 
Hence $A_{C_1}-A_{C_2}=\sum_r b_r\sqrt{\Delta}E_r$. \qed

\begin{remark}
Note that, unlike the abelian group case, the feasible $\Delta$s are not determined by the group exponent: the exponent only provides all the candidates for $\Delta$, some of which might not be achievable.  That is because for some $\Delta$, 
we have $K_j^k=K_j$ for all conjugacy classes $K_j$ and all $k\in \Z_{n_1}^k$.  In this case if $C_1$ is a union of conjugacy classes, then for $k\in \Z_{n_1}^*\backslash H_\Delta$ and  $C_2=C_1^k$, we in fact have $C_2=C_1$, hence the two sets in part (ii) of the theorem are identical, contradicting the hypothesis that they are disjoint. 
\end{remark}
For example,  both the symmetric group $S_3$ and the alternating group $A_4$ have group exponents equal to 6. There exists an oriented normal Cayley graph for  $A_4$, all of whose eigenvalues are integer multiple of $\sqrt{\Delta}$ (see the example below), while such an oriented normal Cayley does not exist on $S_3$.  
In fact, there is no oriented normal Cayley graph for the symmetric group $S_n$, as each conjugacy class is closed under inverse. 
\begin{example}
Consider the alternating group $A_4$. It has exactly four conjugacy classes: $K_0=\{(1)\}, K_1=\{(12)(34),(13)(24),(14)(23)\}, K_2=\{(123),(243),\\ (134),(142)\}$ and $K_3=\{(132),(234),(143),(124)\}=K_2^{-1}$. 
The group exponent of $A_4$ is 6. Therefore the only possible $\Delta$ is $-3$. 
Direct computation shows that $H_{-3}=\{1\}\leq \Z_6^*$ and 
the oriented Cayley graph $X(A_4, K_2)$ (or $X(A_4, K_3)$) has eigenvalues $4\sqrt{-3}, -4\sqrt{-3}$ (simple) and 0 of multiplicity 10. 
\end{example}
\begin{example}
The alternating group $A_5$ has exactly five conjugacy classes, with respective representatives $(1)$, $(12)(34)$, $(123)$, $(12345)$ and $(12354)$. 
Denote these conjugacy classes by $K_j$ for $=0,\ldots, 4$.  
Note that $K_0, \, K_1$, and $K_2$ are also conjugacy classes of $S_5$, 
and $K_3\cup K_4$ is the conjugacy class of $S_5$ corresponding the 5-cycles. \\ 
The group exponent of $A_5$ is 30, so the possible $\Delta$s are $1, -3, 5$, and $-15$. 
Since all the conjugacy classes are closed under inverse, 
$\Delta=-3$ and $\Delta=-15$ are not achievable. 
Now $\Gal(\rats(\zeta_{30})/\rats) \cong \Z_{30}^*=\{1, 7, 11, 13, 17, 19, 23, 29\}$ and  direct computation shows  $H_5=\{1, 11, 19, 29\}$. 
The $\sim_1$ equivalence classes are $K_0, K_1, K_2$ and $K_3\cup K_4$. 
When $\Delta=5$, for $j=1,2$, $K_j^k=K_j$ for any $k\in \Z_{30}^*$, hence none of them could be a subset of $C_1$ (otherwise $C_1$ and $C_2=C_1^k$ intersect nontrivially). 
Direct computation shows that $K_3$ and $K_4$ are both unions of $\sim_5$ equivalence classes, and that $K_3^k=K_4$ for  $k\in \Z_{n_1}\backslash H_5$. 
Hence all eigenvalues of the signed Cayley $X(A_5, K_3, K_4)$ are  integer multiples of $\sqrt{5}$.\\
The semidirect product group $G=C_7\rtimes C_3$ has five conjugacy classes, corresponding to elements of order 1, 3, 3, 7, and 7, respectively. 
For $k=3,7$, by use of a conjugacy class corresponding to elements of order $k$, one obtains an oriented normal Cayley graph for $G$ all of whose eigenvalues are integer multiple of $\sqrt{-k}$; by taking the union of the two conjugacy classes  corresponding to element of order $k$, one get integral Cayley graphs.   There are no signed Cayley graph for $G$ all of whose eigenvalues are integer multiple of $\sqrt{21}$. 

\end{example}
\begin{remark}
Recall that $n_1$ is the exponent of the group $G$ and $\cL=\rats(\zeta_{n_1})$. 
Let $L$ be the field extension of $\rats$ generated by the irreducible character values of $G$ (or equivalently, the eigenvalues of the conjugacy class scheme of $G$). 
Then $L\subseteq \cL$ and $L$ determines the exact set of feasible $\Delta$s:  
these are the $\Delta$s such that $\rats(\sqrt{\Delta})\subseteq L$,  
corresponding to the subgroups $\Gal(L/\rats(\sqrt{\Delta}))$ of $\Gal(L/\rats)$ of index 2.  
For abelian groups, the two fields $L$ and $\cL$ coincide. 
But if $L$ is not a cyclotomic field, it might be harder to find the $\sim_\Delta$ equivalence relation on $G$ directly, so we consider the larger field $\cL$ and make use of the group $\Gal(\cL/ \rats(\sqrt{\Delta}))\cong H_\Delta  \subseteq \Z_{n_1}^*$. 
\end{remark}

Note that a part of the result in Section~\ref{sec:mainthm} on abelian groups can be deduced from Theorem~\ref{thm:charnonab}, as all the conjugacy classes of an abelian group have size 1. The proof presented in Section~\ref{sec:mainthm} does not require explicit knowledge of characters or  representations, by working with the exact values of the idempotents $E_i$s of the  scheme.  Also note that the result is about normal Cayley graphs, not about general Cayley graphs on nonabelian groups.

We finish this section by providing examples of an oriented Cayley graph and a signed Cayley graph for $S_4$, both connected non-normal, whose eigenvalues are all integer multiple of $\sqrt{\Delta}$. They also serve as examples to show that for a nonabelian group $G$ with exponent $n_1$, there might exist an oriented or signed Cayley graph for $G$ such that all its eigenvalues are integer multiple of $\sqrt{\Delta}$ with $\rats(\sqrt{\Delta})\not\subseteq \rats(\zeta_{n_1})$. This differs from normal Cayley graphs. 
\begin{example}
Consider the symmetric group $S_4$. Its exponent is equal to 12, so all the possible $\Delta$s for a normal Cayley graph for $S_4$ are $\pm1, \pm3$. 
The subgroup generated by $C=\{(1234), (1243)\}$ is the whole group $S_4$. 
All the eigenvalues of the oriented Cayley graph $X(S, C)$ are integer multiple of $\sqrt{-2}$ (not in the list). 
One can also check that for $C'=\{(1234),(1243),(1324)\}$, all the eigenvalues of the (connected) oriented Cayley graph $X(S_4, C')$ are integer multiple of $\sqrt{-3}$; 
for $C_1=\{(1234), (1432),(1243),(1342)\}$ and $C_2=\{(1324),(1423)\}$, the (connected) signed Cayley graph $X(S_4,C_1,C_2 )$ has only integer eigenvalues; 
all the eigenvalues of the (disconnected) signed Cayley graph $X(S_4,\{(123),(132)\},\{(124),(142)\})$ are integer multiple 
of $\sqrt{2}$. 
\end{example}

\section{Continuous quantum walks}
\label{sec:umpstexamples}

In the previous sections, we have found a characterization of the possible 
connection sets for an oriented (signed) Cayley graph to have  
eigenvalues in a quadratic extension of the rationals. 
In this section, we apply our knowledge to continuous 
time quantum walks. 

Let $M$ be a Hermitian matrix associated to a weighted graph $X$. 
The transition matrix of the continuous time quantum walk on $X$ relative to $M$ at time $t$ is 
\[
	U_X(t)=\exp(itM). 
\]
The walk on $X$ is \textsl{periodic} at vertex $a$ at time $t$ if 
$U_X(t)e_a=\alpha e_a$ for some complex number $\alpha$. 
We say $X$ admits \textsl{perfect state transfer} from 
vertex $a$ to $b$ at time $t$ if 
$U_M(t)e_a=\alpha e_b$, 
admits \textsl{local uniform mixing} relative to vertex $a$ if $U(t)e_a$ 
is flat, where a vector or a matrix is said to be flat if all its entires have the same modulus, 
and admits \textsl{uniform mixing} at time $t$ if $U(t)$ is flat. 
An eigenvalue $\theta$ is said to be in the \textsl{eigenvalue support} 
of vertex $a$ if there exists an eigenvector $y$ of $M$ associated to $\theta$ such that $e_a^Ty\neq 0$. 
Sometimes we refer to the walk on $X$ relative to $M$ as the walk on $M$ for short. 

The phenomenon of periodicity is in fact closely related to the 
spectrum of $X$: 
the continuous time quantum walk on an oriented graph $X$ is periodic at vertex $a$ if and only if 
all the eigenvalues in the support of $a$ are  integer multiples of $\sqrt{\Delta}$ 
for some square-free integer $\Delta$ \cite{GLori}. 
For (oriented, signed, or unsigned) Cayley graphs, 
periodicity at a vertex in fact requires all eigenvalues of the graph 
to be integer multiples of $\sqrt{\Delta}$ (since all the eigenvalues of such a 
graph are in the eigenvalue support of every vertex). 
Here we prove a more general result on roots-of-unity-weighted walk-regular graphs. 
The proof uses ideas in Section 3 of \cite{godsil2011peri};
we include a proof for completeness. 
Let $X$ be a weighted graph on $n$ vertices with adjacency matrix $A$.  
We say $X$ is \textsl{weighted-walk-regular} if $A^k$ has constant diagonal for $k=1,\ldots, n-1$. 
That is, 
if all the $n$ principal minors of $tI-A$ of size $(n-1)\times (n-1)$ are equal. 
When talking about quantum walks on a weighted graph $X$, 
we use a Hermitian matrix $M$ associated to the graph. 

\begin{lemma} \label{lem:peridel} 
Let $X$ be a weighted-walk-regular graph with all weights $m$-th 
roots of unity for some integer $m$.  
Then the continuous quantum walk on $X$ is periodic at a vertex 
(and therefore at all vertices) if and only if all eigenvalues of $X$ 
are integer multiple of $\sqrt{\Delta}$ for some square-free integer $\Delta$. 
\end{lemma}

\proof
Assume $X$ has $n$ vertices, 
and denote the Hermitian matrix associated to $X$ as $M$. 
Let $\{\theta_1,\ldots, \theta_{n}\}$ be the multi-set of eigenvalues of $M$. 
Note that they are all algebraic integers, 
as eigenvalues of a matrix all of whose entries are algebraic integers are 
themselves algebraic integers. 

If $X$ is periodic at vertex $a$ at time $t$, 
then since $X$ is walk-regular, $U_X(t)=\gamma I$ for some complex number 
$\gamma$ with $|\gamma|=1$ and $X$ is in fact periodic 
at all vertices at time $t$. 
Since $M$ has zero-diagonal ($X$ has no loops), 
the sum of all eigenvalues of $M$ is 0. 
Therefore $\gamma^{n}=\det\left(U_X(t)\right)=\prod_{j=1}^{n}e^{it\theta_j}=e^{it\tr(M)}=1$ 
and $U(nt)=(\gamma I)^{n}=I$.  
Hence $e^{int\theta_j}=1$ for any $j=1,\ldots, n$ 
and $nt\theta_j=2m_j\pi $ for some integers $m_j$.  
Therefore for any two eigenvalues $\theta_r, \theta_s$ of $M$ with $\theta_s\neq 0$, 
\[
\frac{\theta_r}{\theta_s}\in \rats. 
\]
Assume $\theta_1\neq0$, 
then there exist rationals $q_r$ such that $\theta_r=q_r\theta_1$. 
Since $M$ is Hermitian with all its entries $m$-th roots of unity, 
$\sum_{j=1}^{n}\theta_j^2=\tr(M^2)$ is twice the number of edges of $X$. 
Hence 
$(\sum_r q_r^2)\theta_1^2=\sum_r\theta_r^2 $
is an integer. 
Therefore $\theta_1^2$ is a rational algebraic integer and therefore an integer. 
Assume $\theta_1=m \sqrt{\Delta}$ for some integer $m$ and some square-free integer $\Delta$. Then the fact $\theta_r=q_r\theta_1=q_rm\sqrt{\Delta}$ is an algebraic integer implies 
that $q_rm\in \Z$ and therefore all eigenvalues of $M$ are integer multiples of $\sqrt{\Delta}$. \qed

From the above result we know that if all weights of a weighted Cayley graph $X$ are $m$-th roots of unity, 
then $X$ is periodic (at a vertex) if and only if all its eigenvalues are 
integer multiple of $\sqrt{\Delta}$. 
Now we characterize when a $4$-th roots of unity weighted Cayley graph $X$ on 
an abelian group is periodic. 
The graph is in fact a combination of an oriented Cayley graph $X(G,C(i))$ 
and a signed Cayley graph $X(G, C(1), C(-1))$, with $C(i)$, $-C(i)$, $C(1)$ 
and $C(-1)$ mutually disjoint. 
The underlying graph is $X(G, C(i)\cup (-C(i))\cup C(1) \cup C(-1))$.   
The arcs corresponding to $C(i)$ receive weights $ i$, 
the edges corresponding to $C(1)$ receive weights $1$, 
and the edges corresponding to $C(-1)$ receive weights $-1$.  
We denote such a graph as $X(G,C(i),C(1),C(-1))$, 
and call it a \textsl{mixed} graph. 

\begin{theorem}\label{thm:periherm}
Let $G$ be a finite group with exponent $n_1$, 
and let $\Delta$, $H_\Delta$ and the equivalence relations $\sim_\Delta$ 
be defined as in Theorem~\ref{thm:charnonab}. 
Let $X=X(G,C(i),C(1),C(-1))$ be a weighted normal Cayley graph for $G$ ($C(i), C(i)^{-1}, C(1),C(-1)$ are unions of conjugacy classes and pairwise disjoint) with weights $i$, 
$-i$, $1$ and $-1$ corresponding to $C(i)$, $C(i)^{-1}$,  $C(1)$, and $C(-1)$, 
respectively. Then the continuous quantum walk on $X$ is periodic if and only if 
for some $\Delta$ in \eqref{eq:deltas}, 
all of the following holds :
\begin{itemize}
\item 
$C(i)$ is a union of $\sim_{-|\Delta|}$ equivalence classes, 
\item 
$C(1)$ and $C(-1)$ are unions of $\sim_{|\Delta|}$ equivalence classes, 
\item 
and further $C(1)=C(-1)^k$ for some (any) $k\in \Z_{n_1}^*\backslash H_{|\Delta|}$ if $\Delta \neq 1$. 
\end{itemize}
We adopt the convention $C(i)=\emptyset$ if $-|\Delta|$ does not 
satisfy \eqref{eq:deltas} (in this case, $X$ is a signed Cayley),  
$C(1)=C(-1)=\emptyset$ if $|\Delta|$ does not satisfy \eqref{eq:deltas} 
($X$ is an oriented Cayley graph).  
\end{theorem}

\proof 
Since $X$ is weighted-walk-regular (in fact weighted-vertex-transitive),  
by Lemma~\ref{lem:peridel}, 
$X$ is periodic if and only if 
all its eigenvalues are integer multiples of $\sqrt{\Delta}$ for some 
square-free integer $\Delta>0$ (as eigenvalues of a Hermitian matrix are real). 
This requires all eigenvalues of the normal oriented Cayley $X(G, C(i))$ have the 
form $m_r\sqrt{-\Delta}$, and all eigenvalues of the normal signed 
Cayley graph $X(G, C(1),C(-1))$ has the form $n_r\sqrt{\Delta}$. 
Take the corresponding set to be empty if $-\Delta$ or $\Delta$ does not satisfy \eqref{eq:deltas}. 
Now the result follows from Theorem~\ref{thm:charnonab} and Lemma~\ref{lem:peridel}.\qed

It is known that the continuous quantum walk on a Hermitian matrix $M$ is 
periodic at vertex $a$ if and only if the following ratio condition on the eigenvalues of $M$ 
holds \cites{URAcomp,godsil2017real}: 
for any $\theta_r,\theta_s,\theta_k,\theta_\ell$ in the eigenvalues 
support of $a$ with $\theta_k \neq \theta_\ell$, 
\[
\frac{\theta_r-\theta_s}{\theta_k-\theta_\ell}\in \rats.
\]
Further, Godsil shows that if the continuous quantum walk on a Hermitian 
matrix with all entries algebraic or all entries real admits perfect state 
transfer from vertex $a$ to $b$, then the above ratio condition holds 
on the eigenvalue support of $a$ \cite{godsil2017real}. 
Therefore given a Hermitian matrix with algebraic entries or given a real 
symmetric matrix, periodicity at vertex $a$ is a necessary condition for 
perfect state transfer from $a$ to some other vertex to occur. 
Lemma~\ref{lem:peridel} characterizes when a root-of-unity-weighted walk-regular 
graph is periodic, and Theorem~\ref{thm:periherm} characterizes such  
mixed Cayley graphs in terms of the connection sets. 
Furthermore, uniform mixing 
 relative to $a$ in an oriented graph implies that $X$ is periodic at vertex $a$ \cite{godsil2017real}. 
This tells us that for interesting quantum walk phenomenon to occur on an 
oriented Cayley graph $X$, all eigenvalues of $X$ are of the 
form $m\sqrt{\Delta}$ for a fixed square-free integer $\Delta$. 
Our characterization in Theorem~\ref{thm:char} provides us the graphs to look at for interesting quantum walk phenomenon to occur. 
In the following we give some oriented Cayley graphs that admit uniform mixing or multiple state transfer.

\subsection{Uniform mixing}

Let $X$ be an oriented Cayley graph on $n$ vertices 
and denote its skew-symmetric adjacency matrix by $A$. 
Then the transition matrix of the continuous time quantum walk on $X$ at time $t$ is 
\[ 
U_X(t)=\exp(it(iA))=\exp(-tA),
\]
a real orthogonal matrix. 
If $X$ admits uniform mixing at time $t$, 
then all the entries of the unitary matrix $U(t)$ are equal to $\pm\frac{1}{\sqrt{n}}$.  
Therefore $\sqrt{n}U(t)$ is a real Hadamard matrix, 
and hence $n=1,2$ or  $4k$ for some positive integer $k$. 
Furthermore,  
if $X$ is an oriented circulant (the associated abelian group $G$ is cyclic)  
admitting uniform mixing at time $t$, 
then $U(t)$ is a circulant Hadamard matrix, 
which is conjectured to exist only for $n=1, 4$ \cite{Ryser}. 

In general, we have

\begin{lemma}\label{lem:umevensq}
	Let $X$ be an oriented Cayley graph for a group of size $n>2$.
	If $X$ admits admits uniform mixing, 
	then $n$ is an even perfect square. 
\end{lemma}

\proof
If $X$ admits uniform mixing at time $t$, then $U(t)$ is flat.
Since it is real, $H=\sqrt{n}U(t)$ is a real Hadamard matrix and,
as it lies in the group algebra, it is regular 
(all row and column sums are equal). Therefore
there is a constant $c$ such that $H\one=c\one$ 
and $n\one=H^TH\one=c^2\one$. It follows that $n$ is a perfect square 
and, as $H$ is Hadamard matrix of order greater than two, its order 
is divisible by $4$.\qed

Note that the above result holds for any oriented Cayley graph, including non-normal Cayley graphs. 

\begin{corollary}
There is no oriented Cayley graph for the alternating group $A_n$ or the symmetric group $S_n$ on which uniform mixing occurs. 
\end{corollary}
\proof 
The result follows from the fact that the group sizes  $n!/2$ (of $A_n$) and $n!$ (of $S_n$) are not perfect squares and Lemma~\ref{lem:umevensq}. \qed

Let $X_1=X(\Z_4,\{1\})$ be the oriented cycle on 4 vertices. 
Direct computation shows that $U_{X_1}(\pi)=I_4$,
\[
U_{X_1}(\frac{\pi}{4})=\frac12
\begin{bmatrix}
1 & 1 & 1 & -1\\
-1 &1 & 1 & 1\\
1 & -1 &1 & 1\\
1 & 1 & -1 & 1
\end{bmatrix},
U_{X_1}(\frac{\pi}{2})=
\begin{bmatrix}
0 & 0 & 1 & 0\\
0 &0 & 0 & 1\\
1 & 0 &0 & 0\\
0 & 1 & 0 & 0
\end{bmatrix}. 
\]
Therefore  
the directed cycle on 4 vertices admits uniform mixing at 
time $\frac{\pi}{4}$, 
and admits perfect state transfer at time $\frac{\pi}{2}$. 
As $U(\frac{3\pi}{4})=U(\frac{\pi}{4})U(\frac{\pi}{2})$, 
$X_1$ also admits uniform mixing at time $\frac{3\pi}{4}$.

Now let $X=X(\Z_4^n,\{(1,0,\ldots, 0),(0,1,0,\ldots, 0),\ldots, (
0,\ldots, 0,1)\})\cong X_1^{\square n}$. 
Then 
\[
A(X)=A(X_1)\otimes I_4\otimes\cdots\otimes I_4+ 
     I_4\otimes A(X_1)\otimes I_4\otimes\cdots\otimes I_4+\cdots
     I_4\otimes\cdots\otimes I_4\otimes A(X_1). 
 \]
At any time $t$, 
\[
U_X(t)=U_{X_1}(t)\otimes U_{X_1}(t)\otimes \cdots\otimes U_{X_1}(t).
\]
In particular,
$U_X(\frac{\pi}{4})=U_{X_1}(\frac{\pi}{4})^{\otimes n}$ is flat,
that is, the oriented Cayley $X$ admits uniform mixing at time $\frac{\pi}{4}$. 
In fact, in the above connection set, 
we could also change any number of 1s to $-1$s and preserve uniform mixing. 
These are trivial examples of graphs with uniform mixing obtained from 
the Cartesian product of smaller oriented graphs with uniform mixing. 
It was shown when one orients a bipartite graph according to the bipartition (every arc goes from one part of the bipartition to the other part), 
then the resulting oriented graph has the same quantum information transfer 
property as the original bipartite graph, as for any time $t$, the two 
transition matrices are similar \cite{godsil2017real}. 
Here we give some more interesting examples. 

\begin{example}
The oriented Cayley graph 
$X=X(\Z_4^2,\{(1,0),(0,1),(1,1)\})$ admits uniform mixing at time $\frac{\pi}{4}$. 
Note that $X$ is not a Cartesian product, 
and its underlying graph is not bipartite. 
More generally, 
\[
X=X(\Z_4^{2k},\{1,0,\ldots, 0),(0,1,0,\ldots, 0),\ldots, (0,\ldots, 0,1),(1,\ldots,1\})
\]
admits uniform mixing at time $\frac{\pi}{4}$, 
but is neither a Cartesian product nor bipartite. 
\end{example}

\subsection{Multiple state transfer}

Let $X$ be an oriented graph. 
If there is a subset $S$ of $V(X)$ with $|S|>2$ 
such that between any two vertices of $S$ there is perfect state transfer, 
then we say $X$ admits \textsl{multiple state transfer} on $S$ \cite{GLori}. 
The only known examples of multiple state transfer are the oriented 3-cycle and 
an oriented graph on 8 vertices \cite{GLori} and \cite{URAcomp}, 
and the ones constructed in \cite{URAcomp}. 
Here making use of our characterization of oriented Cayley graphs with nice eigenvalues, 
we provide several infinite families of oriented Cayley graphs with multiple state transfer. 

\begin{example}
We show some examples of oriented circulant with multiple state transfer here. 
\begin{itemize}
\item
Let $m\geq 3$ and $C=\{2^{m-d}(4r+1) \; |\,  d=2,\ldots, m, \;  r=0,1,\ldots, 2^{d-2}-1\}$, then the oriented Cayley graph (cocktail party graph) $X(\Z_{2^m},C)$ 
admits multiple state transfer on $S=\{0,2^{m-2}, 2\times 2^{m-2}, 3\times 2^{m-2}\}$ 
and the shifted sets $a+S$ for $a\in \Z_{2^m}$. 
In particular, $X(\Z_8,\{1,2,5\})$ is switching equivalent to the 
example of multiple state transfer in \cite{GLori}. 
\item 
Let $C= \{3^{m-d}r \; |\,  d=1,\ldots, m, \, r=0,1,\ldots, 3^{d-1}-1\}$, 
then the oriented Cayley graph (tournament) 
admit multiple state transfer on $\{0, 3^{m-1},2\times 3^{m-1}\}$ and all the shifted sets.

\item
Let $a$ be any positive integer, 
and $C=\{2^{m-3}a,2^{m-3}2a,2^{m-3}5a\}\cup
\{b\in \Z_{2^ma}^*\; |\, b=4k+1\}$, 
then the oriented Cayley graph 
$X(\Z_{2^ma}, C)$ 
admits multiple state transfer on 
$\{0, 2^{m-2}\times a,2\times 2^{m-2}\times a,3\times 2^{m-2}\times a\}$
 and their shifts.
\end{itemize}
\end{example}

\section{Further questions}
Cayley graphs (oriented or signed) with all eigenvalues integer multiples of $\sqrt{\Delta}$ are the Cayley graphs on which the continuous quantum walk is periodic, which is a necessary condition for perfect state transfer to occur. 
In this paper, we obtained a characterization of when an oriented or signed normal Cayley graph has such nice eigenvalues, in terms of the connection set. We provided examples of oriented or signed non-normal Cayley graphs with such eigenvalues. It would be interesting to get a characterization for all the Cayley graphs. Our proof by use of association schemes for normal Cayley graphs does not apply in this case. 

\begin{question}
Is there a simple characterization of the connection set of an oriented or signed non-normal Cayley graph all of whose eigenvalues are integer multiple of $\sqrt{\Delta}$ for some square-free integer $\Delta$? 
\end{question}

As an application of our characterization, we found infinitely many oriented Cayley graphs which admit multiple state transfer. All the examples we have so far for multiple state transfer are on subsets of size three or four. 

\begin{question}
Is there an unweighted oriented graph $X$ which admits multiple state transfer on some subset $S\subseteq V(X)$ with $|S|>4$?
\end{question}

\end{document}